\documentclass{article}

\usepackage{amsmath}
\usepackage{amssymb}
\usepackage{amsthm}
\usepackage{amsrefs}
\usepackage{hyperref}

\theoremstyle{plain}%
\newtheorem{theorem}{Theorem}
\newtheorem{lemma}{Lemma}
\newtheorem{corollary}{Corollary}

\theoremstyle{remark}%
\newtheorem{remark}{Remark}%

\theoremstyle{definition}%
\newtheorem{definition}{Definition}%

\newcommand{\s}{\sigma}
\renewcommand{\eqref}[1]{(\textup{\ref{#1}})} 

\newcommand{\inn}[2]{\left\langle#1,\,#2\right\rangle}
\DeclareMathOperator{\ran}{ran}

\DeclareMathOperator{\dist}{dist}

\newcommand{\grad}{\nabla}
\newcommand{\Lap}{\Delta}
\newcommand{\di}{\partial}

\newcommand{\si}{\sigma}
\newcommand{\eps}{\epsilon}
\newcommand{\ls}{\lesssim}
\newcommand{\g}{\gamma}
\newcommand{\al}{\alpha}
\newcommand{\Si}{\Sigma}

\renewcommand{\b}{\bar}
\newcommand{\Cb}{\mathbb{C}}

\newcommand{\Rb}{\mathbb{R}}

\newcommand{\E}{\mathcal{E}}
\newcommand{\J}{\mathcal{J}}

\newcommand{\om}{\omega}

\renewcommand{\b}{\bar} 

\newcommand{\abs}[1]{\left\lvert#1\right\rvert}
\newcommand{\norm}[1]{\left\lVert#1\right\rVert}

\newcommand{\Set}[1]{\left\{#1\right\}}
\newcommand{\fullfunction}[5]{\ensuremath{
		\begin{array}{ccrcl}
			{#1}    & \colon  & {#2} & \longrightarrow & {#3} \\
			\mbox{} & \mbox{} & {#4} & \longmapsto     & {#5}
\end{array}}}
\newcommand{\pd}[2]{\ensuremath{\frac{\partial#1}{\partial#2}}}
\newcommand{\od}[2]{\ensuremath{\frac{d#1}{d#2}}}
\newcommand{\md}[6]{\ensuremath{
		\ifinner
		\frac{\partial{^{#2}}#1}{\partial{#3^{#4}}\partial{#5^{#6}}}
		\else
		\frac{\partial{^{#2}}#1}{\partial{#3^{#4}}\partial{#5^{#6}}}
		\fi
}}
\newcommand{\del}[1]{\left(#1\right)}
\newcommand{\thmref}[1]{Theorem~\ref{#1}}

\newcommand{\defnref}[1]{Definition~\ref{#1}}
\newcommand{\secref}[1]{Section~\ref{#1}}
\newcommand{\lemref}[1]{Lemma~\ref{#1}}

\newcommand{\remref}[1]{Remark~\ref{#1}}

\begin{document}
	
		
		
		
		\title{A Generic Framework of Adiabatic Approximation for
			Nonlinear Evolutions}
		
		\author{Jingxuan Zhang}
		
		\maketitle
		
		\begin{abstract}
			In the study of evolution equations, the method of adiabatic approximation is an essential tool to reduce an infinite-dimensional dynamical system to a simpler, possibly finite-dimensional one. In this paper, we formulate a generic scheme of adiabatic approximation that is valid for an abstract nonlinear evolution under mild regularity assumptions. The key prerequisite for the scheme is the existence of what we call approximate solitons. These are some low energy but not necessarily stationary configurations. The approximate solitons are characterized by a number of parameters (possibly infinitely many), and have a manifold structure. The adiabatic scheme reduces the given abstract evolution equation to an effective equation on the manifold of approximate solitons. We give sufficient conditions for the approximate solitons so that the reduction scheme is valid up to a large time. The validity is determined by the energy property of the original evolution as well as the adiabaticity of the approximate solitons. 
		\end{abstract}

		\section{Introduction}
		Consider an abstract evolution equation
		\begin{equation}
			\label{1}
			\di_tu=J E'(u).
		\end{equation}
		Here $u=u_t\in U,\,t\ge 0$ is a $C^1$ path of vectors in some open set $U$ in a real Hilbert
		space $X$. 
		The map $E:U\subset X\to \Rb$ is some energy functional which is $C^2$ on $U$. 
		The vector $E'(u)\in X$ is the $X$-gradient of $E$ at $u$.
		
		We assume the operator $J:X\to X$ in \eqref{1} is a bounded
		invertible linear operator, satisfying 
		$$\text{either }J=-1,\quad \text{or }J^*=J^{-1}=-J.$$
		In the first case, $J$ is the negative of 
		the identical map. In the second case, $J$ is a symplectic
		operator. This symplectic condition holds, 
		for example, if $J$  can be represented by the standard symplectic matrix
		$$
		\begin{pmatrix}0&1\\-1&0\end{pmatrix}.
		$$
		These two cases respectively turn \eqref{1} into either a first-order energy dissipative
		dynamics or a Hamiltonian system.
		
		We assume the following global well-posedness result for \eqref{1}:
		\begin{quote}
			For every $u_0\in U$, there exists a path $u\in C^1(\Rb_{\ge0},U)$
			s.th. $u\vert_{t=0}=u_0$ and $\di_tu_t=JE'(u_t)$ for $t>0$.
		\end{quote}
		In this case we say $u_t$ is the flow generated by $u_0$ under \eqref{1}.
		Then we consider the following problem: 
		\begin{quote}
			Suppose a configuration $u_0\in X$ can be parametrized up to a small error
			by a point $\sigma_0$ in a manifold $\Sigma$.
			Let $u_t,\,t\ge0$ be the flow generated by $u_0$ under \eqref{1}. 
			Can one reduce the full flow $u_t\in U$ to
			an approximate flow $\s_t\in\Sigma$ generated by $\s_0$
			under some suitable effective dynamics?
		\end{quote}
		
		Here are some examples when this problem arises:
		\begin{enumerate}
			\item Let $n,\,k\ge1$.
			Let $X$ be a suitable space of functions
			from $\Omega\to \Rb^{k}$, 
			where $\Omega\subset \Rb^{n+k}$ is a domain. 
			Suppose the initial configuration $u_0$ is \textit{localized}
			in the sense that outside a neighbourhood 
			of some $k$-dimensional concentration set
			$\s_0\subset\Omega$,
			all derivatives of $u_0$ vanish rapidly. Then one is interested in whether
			the evolution $u_t$ remains localized near some concentration set
			$\s_t$ for $t>0$, and if so, what kind 
			of geometric flow governs the motion of $\s_t$.
			This problem arises from the study of phase transition, see e.g.
			\cites{MR1100206,MR1177477,MR2032110,MR2195132}.
			
			\item Let $X$ and $u_0$ be as before. Suppose $\s_0$ is given
			by a collection of distinct points in $\Omega$ (i.e. $k=0$). Heuristically, 
			one would expect the evolution
			of $\s_0$ to be the motion law of interaction among the  points
			of localization. One is interested in 
			to what extent the full dynamics can be reduced to this 
			``renormalized'' dynamics of points. 
			This is the setting for various soliton scattering problems in high energy physics \cites{MR1257242,MR1309545,MR2097576,MR2525632}.
			
			\item Let $X$ be a suitable space of geometric objects (e.g. curves, surfaces,
			etc.). Suppose
			$u_0\in X$ is determined by a number of parameters,
			for instance center, radius, axial direction, etc..
			Then one is interested in whether
			$u_t$ can still be described faithfully by these parameters at time $t>0$.
			This problem is essential to the study of rigidity
			under various geometric flows \cites{MR3374960,MR3397388,gang2017dynamics,gang2019mean,zhou2021nondegenerate}.
		\end{enumerate}
		
		\subsection{Outline of the main result} In this study we propose an abstract scheme 
		to answer the questions above, 
		known in the physics literature as the
		method of  \emph{adiabatic approximation}.
		The precise statements of the main results are
		in Theorems \ref{thm1}-\ref{thm2}.
		Below we give an outline.
		
		Our main assumption is the existence of 
		a parametrized family
		of configurations $v_\s\in X,\,\s\in \Sigma$, where $\Sigma$
		is a  manifold. We call $\Sigma$ the \textit{moduli space}, in the sense that $\Sigma$ contains
		various modulation parameters.
		The space $\Sigma$ can be finite-dimensional if, for instance,  it represents finitely many concentration points in a domain. In general, $\Sigma$ can be infinite-dimensional,
		e.g. as  a space of geometric objects or 
		local gauge symmetries.
		We call this family $v_\s$ the \emph{approximate solitons}. 
		In Section \ref{sec:2}, we further explain the terminology, and
		list the precise requirements \eqref{C1}-\eqref{C3} 
		for $v_\s$. These requirements 
		specify certain approximate energy-critical and linear stability properties.
		We give a justification in \secref{sec:2} for the generality of 
		these requirements.
		
		Suppose there exists a family of approximate solitons 
		$v_\s\in X,\,\s\in \Sigma$
		satisfying the main assumptions \eqref{C1}-\eqref{C3}.
		Then the adiabatic approximation scheme goes  
		as follows: First, we find an  evolution equation
		$\di_t\s=F(\s)$ for a
		path $\s_t\in\Sigma$ with the following property:
		Let $u_t,\,t\ge0$ be a solution to \eqref{1}.
		\textit{Assuming $u_t$ remains close} to
		a path $v_\s\equiv v_{\s_t}$,
		then $\s_t$ evolves according to this equation in the leading order.
		We calculate the velocity $F(\s)$ explicitly (see \eqref{eff}).
		This equation for $\s$ is the effective (or adiabatic) dynamics for \eqref{1}.
		In this step we use the approximate critical point
		property \eqref{C1} of  $v_\s$.
		
		At this point there is  a caveat: 
		In the first step, we have assumed that  $u_t$ remains close to
		some approximate soliton for all $t$. But,
		a priori, even $u_0$ lies in a neighbourhood of the approximate solitons,
		the path $u_t$ may soon exit
		this neighbourhood (for instance, due to acceleration). 
		Thus we need to justify
		the validity of the 
		effective dynamics  derived in the first step.
		
		In the second step, we show that so long as the initial configuration $u_0$
		is close to an approximate soliton $v_{\si_0}$, then the flow $u_t$ generated by $u_0$ stays 
		uniformly close
		to a
		path of approximate solitons $v_{\si_t}$,
		at least up to a large time.
		In this step we use  the 
		stability properties \eqref{C2}-\eqref{C3}. 
		This step is analogous to proving
		orbital stability for ground states \cites{MR783974,MR820338}.

		As a corollary of the main result, we derive a converse that 
		allows us to find a flow evolving by \eqref{1}
		that agrees with a given adiabatic flow in $\Sigma$ in the leading order
		up to a large time. 
		The precise statements are given in 
		Corollaries \ref{cor1}-\ref{cor2}.

		\subsection{Historical Remarks}\label{sec:1.2} The idea in this paper dates back
		at least to Manton's moduli space approximation scheme for monopole dynamics \cite{MR647883}. In addition to the cited works above,
		rigorous results using similar methods include 
		\cites{MR1644389,
			MR1644393,MR2094474,MR2232367,
			MR2189216,MR2465296,MR2576380,MR2805465,MR2813581,MR2952171,
			MR3803553,MR4135093,lasser_lubich_2020}, 
		with diverse applications to superfluidity, superconductivity, particle physics, and geometric flows. The point in common in these applications is that 
		solitonic configurations (i.e. coherent states)  arise naturally
		due to focusing nonlinearity or other types of constraints.
		The structure of the present paper can be compared to \cites{MR2094474,MR2232367}.
		
		We refer the readers to an excellent review on the applications of
		adiabatic approximation to classical field theory \cite{MR2360179}.
		The review contains several examples of effective dynamics describing
		interacting point solitons, and discusses a different approach through compactness arguments.
		The latter finds applications to a large class of problems,
		which, among others, include the now-classical geometric theory of phase transitions
		\cites{MR866718, MR1100206,MR1177477, MR2032110}, relating the flow of a 
		real order parameter under the Allen-Cahn equation to the mean curvature flow
		of its nodal set. 
		
		Mathematically, one can compare the results in this paper to
		the classical
		invariant manifold theory
		developed in  \cites{MR635782,MR1000974,MR2439610,MR1445489,MR1675237,MR1160925}.
		In this regard, our main assumptions \eqref{C1}-\eqref{C3} can be compared
		to the normal hyperbolicity condition in those references.
		However, we note that here our focus is different,
		since we are less interested in the properties of the invariant manifold 
		\textit{per se}, but rather to
		reduce \eqref{1} to an explicit effective dynamics
		and make sure the reduction is both tractable and valid at least for a long time.
		
		\subsection{Arrangement} The arrangement of this paper is as follows. In \secref{sec:2},
		we describe the manifold $M$ of approximate solitons.
		We list two groups of assumptions on this manifold
		and discuss the generality of these assumptions for
		applications.
		
		In \secref{sec:3}, we prove two key lemmas.
		By these lemmas, we find a heuristic 
		candidate \eqref{eff} as the effective dynamics of \eqref{1}. 
		
		In Sections \ref{sec:G}-\ref{sec:H},  we
		prove the main results of this paper.
		We show that for both dissipative and 
		Hamiltonian \eqref{1}, the heuristic effective dynamics 
		\eqref{eff} is indeed valid 
		as the adiabatic approximation for \eqref{1} for a long time (global in the dissipative case).
		
		In \secref{sec:5}, we give some concrete example to
		illustrate the application 
		of our abstract framework to the study of the motion of mesoscopic interfaces,
		a central problem in statistical physics.
		In Appendix, we list the basic concepts of variational
		calculus that are used repeatedly.
		
		\subsection*{Notations} Throughout this paper, when no confusion arises, we shall drop the time dependence $t$ in subscripts. 
		An estimate $A\ls B$ means there is some $C>0$ independent 
		of time and all the parameters in question, s.th. $A\le CB$. 
		The expression $A\sim B$ means that
		$A\ls B$ and $B\ls A$ hold simultaneously. 
		
		\section{The Manifold of Approximate Solitons}\label{sec:2}

		The central object of this paper is the manifold $M$ defined below.

		\begin{definition}[manifold of approximate solitons]\label{M}
			Let $X$ be a real Hilbert space.	Let $\Sigma$ be a closed Riemannian manifold sitting in some 
			(possibly
			infinite-dimensional) Banach space.
			Let $f:\Sigma\to U\subset X$ be a $C^2$ map,
			where $U\subset X$ is an open set on which the energy functional $E$ in \eqref{1} is $C^2$.

			Then the  subset $$M:=f(\Sigma)$$ forms a manifold
			in $X$. The tangent space
			$T_{f(\s)}M,\,\si\in \Si$ can be trivialized as a subspace of $X$, as
			$$T_{f(\s)}M=\Set{df(\si)\xi\in X:\xi\in T_\si\Si}.$$
			By assumption, the tangent space $T_\si \Si$ can also be trivialized as a subspace $Y$
			of 
			the ambient Banach space. 
			
			Fix some bases for the tangent spaces $T_{f(\si)}M$ and 
			$T_\si\Si$  w.r.t. the trivializations above.
			Denote by $$g_\s:Y\to X$$  the action of $df(\s):T_\s\Sigma\to T_{f(\s)}M$ on a given fiber.
					In other words, $g_\si$ is the Fr\'echet derivative $df(\si)$ in local coordinates.
			Denote by $$g_\s^*:X\to Y$$
			the adjoint to $g_\s$.
			
			We call  $M$  a \textit{manifold of approximate solitons}
			if the following holds:
			\begin{enumerate}
				\item (Solitonic assumptions) 
				There exist
				$$0<\eps\ll1,\quad \beta>0$$
				such that
				every element  $f(\s)\in M$ with $\si\in \Si$ satisfies the following conditions: 
			\begin{align}
				&\label{C1} \norm{E'(f(\s))}_X\le \eps,\tag{C1}\\
				&\label{C2} \text{$L_\s:=E''(f(\s)):X\to X$ is self-adjoint, and }L_\s\vert_{\del{JT_{f(\s)}M}^\perp}\ge\beta>0,\tag{C2}\\
				&\label{C3} {L_\s}\vert_{T_{f(\s)}M}\le \eps\tag{C3}.
			\end{align}
		In \eqref{C2}, the the bounded invertible operator $J:X\to X$
		is as in \eqref{1}.
		
		\item (Geometric assumptions)
		There exist $$0<c\le C<\infty,\quad 0\le\al<1$$ such that
			the parametrization $f$ satisfies the 
			following conditions: 
			\begin{align}
				\label{5.4'}
				&\text{		$df(\s):T_\si\Si\to X$ is injective at every $\sigma\in \Sigma$,}\tag{G1}\\
				\label{5.4}
				&c\eps^{-\al}\norm{\xi}_Y\le \norm{g_\s\xi}_{Y\to X}\le C\eps^{-\al}\norm{\xi}_Y\quad(\xi\in Y)\tag{G2}.
			\end{align}
			\end{enumerate}
		\end{definition}
		
		\begin{remark}
			Condition \eqref{5.4'} implies that
			$f$ is an immersion of the given manifold  $\Si$.
			Hence, $M$ forms a non-degenerate  manifold in $X$.
		\end{remark}

		Conditions \eqref{C1}-\eqref{5.4} play central roles 
		for the validity
		of adiabatic theory for the full evolution \eqref{1}.
		In the remaining of this section, we discuss these conditions in two groups.
		The first group, \eqref{C1}-\eqref{C3}, concerns with the qualitative properties 
		of the approximate solitons related to the energy functional $E$ in \eqref{1}.
		The second group, \eqref{5.4'}-\eqref{5.4}, concerns with the geometric properties of the manifold 
		$M$ from in \defnref{M}.
		In \secref{sec:5}, we give some examples of approximate solitons
		that arise naturally as models of mesoscopic interfaces.

		\subsection{Solitonic Assumptions}\label{sec:2.1}
		In this subsections, we discuss the solitonic
		assumptions \eqref{C1}-\eqref{C3}.

		Let $M$ be a manifold of approximate solitons
		as in \defnref{M}. Then 
		there exist
		two constants 
		$$0<\eps\ll1,\quad \beta>0$$
		s.th. 
		every element $f(\s)\in M$ with $\si\in \Si$ satisfies 
		\begin{align}
			&\label{C1d} \norm{E'(f(\s))}_X\le \eps,\tag{C1}\\
			&\label{C2d} \text{$L_\s:=E''(f(\s)):X\to X$ is self-adjoint, and }L_\s\vert_{\del{JT_{f(\s)}M}^\perp}\ge\beta>0,\tag{C2}\\
			&\label{C3d} {L_\s}\vert_{T_{f(\s)}M}\le \eps\tag{C3}.
		\end{align}
	
				\emph{Through out the remaining sections, $\eps$ is taken to be
			a sufficiently small parameter.} In practice, this parameter usually comes with
		the equation \eqref{1}.
		
		Condition \eqref{C1} says that elements in $M$ are approximate critical 
		points of $E$. Conditions  \eqref{C2}-\eqref{C3} 
		assert certain stability properties regarding the linearized operator $L_\s$.
		(See  Appendix for the definition of the Hessian $E''$.)
		
		\begin{remark}
			If $E$ is a $C^2$ map on an open set $U\subset X$,
			and $X$ lies in a larger Hilbert space $Z$ with a possibly weaker topology,  
			then the $Z$-gradient $E'(u),\,u\in U$ is an element in $Z$.
			In this situation, the (possibly nonlinear) 
			map $u \mapsto E'(u)$ is $C^1$ from $X\to Z$,
			and for fixed $u\in U$,  the linearized operator $E''(u)$ is bounded from $X\to Z$.
			This would result in a number of changes in the conditions \eqref{C1}-\eqref{C3}.
			For example, we would need estimates on $\norm{E'(f(\si))}_Z$ in \eqref{C1}
			and $\norm{L_\si}_{X\to Z}$ in \eqref{C2}. 
			
			However, in most applications where this situation arises,
			it is possible to show that on the class
			of configurations of interest, the formally weaker topology induced
			by the $Z$-norm is equivalent to that induced by the $X$-norm (e.g. all cited works in \secref{sec:1.2}).
			For example, consider the typical situation where $E'$ is an elliptic 
			differential operator of second order
			with sufficiently regular coefficients (but not necessarily linear)
			mapping from some Sobolev space  $H^k\to H^{k-2},\,k\ge2$. Then, so long as 
			the initial configuration $u_0$
			for \eqref{1} is sufficiently regular, by standard elliptic regularity
			theory we have $\norm{u_t}_{H^k}\ls \norm{u_t}_{H^{k-2}}$
			for the flow $u_t,\,t\ge0$ generated by $u_0$ under \eqref{1}.
			Moreover, the Fr\'echet derivatives depend only on the topology
			but not the norm on the ambient space \cite{MR1336591}*{Sect. 1}.
			Hence, 
			the technicality mentioned
			above has little to no bearing 
			for most applications of  the adiabatic theory developed in the present
			paper, as it can be easily amended \textit{ad hoc}.
			For this reason, in the sequel
			we choose not to pursue this point any further so as not to obscure our main results.
		\end{remark}
		
		Now we explain some terminology.
		We call the space 
		$M$ in \defnref{M} the manifold of \textit{approximate solitons}, and $\Sigma$ the\textit{ moduli space,} for the following reason:
		Consider the case for $J=-1$.
		Suppose $M$ consists 
		of exact critical points of $E$. Then  \eqref{C1}
		holds with $\eps=0$. If, moreover,  $\Sigma$ is the (continuous)
		symmetry group 
		of $E$ and $f$ is the action of $\Sigma$ on $X$,
		then \eqref{C2} holds if the elements in $M$
		are (linearly) stable,
		and \eqref{C3} holds with $\eps=0$, since 
		in this case the tangent space $T_{f(\si)}M$ consists exactly of the
		zero-modes generated by the broken symmetries.

		In the preceding situation, each element in $M$ is an \textit{exact soliton},
		and $\Si$ consists of the \textit{modulation parameters} of the continuous symmetries.
		Conditions \eqref{C1}-\eqref{C3} relax this limiting case,
		while  retaining the key  qualitative features.
		This explains our terminology as well as the motivations for \eqref{C1}-\eqref{C3}. 
		
		We note that the idea of approximating a flow near $M$ by a flow on the
		moduli space dates back to Manton's classical work \cite{MR647883}, and 
		such approximating scheme was first rigorously implemented in \cites{MR1257242,MR1309545}.
		
		\begin{remark}\label{rem2}
			Among the three main conditions above, the most restrictive one 
			is \eqref{C2}. This condition asserts that there is a spectral gap at
			$0$ for the linearized operator $L_\si$. When the operator $L_\si$
			does not have compact resolvent, for instance due to the non-compactness of 
			domain on which $L_\si$ acts,
			verifying this spectral gap condition is a non-trivial issue.
			This affects the applicability of our framework
			to problems arising from e.g. nonlinear optics, where the linearized
			operators at exponentially decaying ground states on $\Rb^d$
			in general possess delicate spectral properties.
			
			Nonetheless,
			we remark that linearly stable exact solitons defined on non-compact
			domains are readily available in various important models.
			For applications to classical field theory, due to the 
			Higgs mechanism, one can obtain coercivity estimates
			of the form \eqref{C2} for the linearized operator at 
			ground states \textit{even in the presence of essential spectrum}. See
			\cites{MR1257242,MR1309545,MR1479248,MR1772246}
			for some examples in gauge field theory of this kind.
			For applications in geometric analysis, due to
			the topological properties of the underlying domain, 
			\eqref{C2} can also be achieved by the linearized operators at key solutions
			defined on non-compact manifolds.
			See \cites{MR2993752,MR3374960} for some examples arising from
			the analysis of the mean curvature flow.
		\end{remark}
		
		Now we proceed to give some further  justifications of the generality of
		the conditions \eqref{C1}-\eqref{C3}.
		Indeed, these conditions 
		are generic, first and foremost, because if one has a single approximate soliton $u_0$
		satisfying these conditions to begin with,
		then one can find other approximate solitons  by perturbing 
		this $u_0$.
		These perturbations can be rather arbitrary,
		because  \eqref{C1}-\eqref{C3} do not impose any quantitative constraints
		other than that the perturbations have to be small.

		For simplicity, suppose $J=-1$ and $u_0=0\in X$ is a given approximate
		soliton, satisfying 
		$$\norm{E'(0)}_X\le \eps,\quad L_0P^+\ge\beta,\quad L_0(1-P^+)\le\eps\quad(L_0:=E''(0)),$$
		where $P^+$ denotes the Riesz projection onto
		the intersection of $\Rb_{>0}$ with the spectrum of $L_0$.
		Consider an arbitrary perturbation $v\in X$ around $u_0=0$
		with $\norm{v}_X\ll1$.
		For \eqref{C1}, the size of
		$\norm{E'(v)}_X$ is still small by the continuity of $E'$.
		For \eqref{C2}, if the energy functional 
		$E$ is sufficiently regular, say at least $C^4$, then the path $v\mapsto L_v:=E''(v)$ is 
		$C^2$, and it follows from \cite{MR2006008}*{Thm. (D)} that
		the second lowest positive eigenvalue of $L_v$
		is at least $\beta/2$ for every $v$ along  this path,
		provided  $\norm{v}_X\ll1$ and 
		the linearized operators $L_v$ have compact resolvents
		(which, of course, excludes the problem due to essential spectrum mentioned in \remref{rem2}). A similar statement holds for  \eqref{C3}.
		
		Hence, if $u_0=0$ is an approximate soliton,
		and the linearized operators have compact resolvents,
		then any element in an immersed manifold $M\subset B_\delta(0)\subset X$ with 
		$0<\delta\ll1$ and 
		tangent space $T_v M$ close to $\ran(1-P^+)$ satisfies \eqref{C1}-\eqref{C3}.
		The condition on the tangent space determines the geometry of $M$ (c.f. the normal
		hyperbolicity condition for infinite-dimensional
		invariant manifolds in \cites{MR1000974,MR2439610,MR1445489,MR1675237,MR1160925}).
		
		Note that  
		even if the linearized operators do not have compact resolvents,
		it is still possible to construct such a manifold, see e.g. \cites{MR1257242,MR1309545,MR2189216}.
		Note also that in the above discussion, $u_0$ needs not to be an exact soliton to begin with,
		see e.g. \cites{MR2097576,MR2189216}.

		To summarize, our emphasis here is the freedom in constructing $M$,
		since we do not require \eqref{C1} with $\eps=0$,
		and \eqref{C2}-\eqref{C3} are persistent under sufficiently regular perturbations
		(in some cases even in the presence of essential spectrum).

		\subsection{Geometric Assumptions}\label{sec:3.2}
		The secondary assumptions are about the parametrization $f$
		of the manifold of approximate solitons given in \defnref{M}.
		
		Recall that we have defined
		$g_\si:Y\to X$ as the Fr\'echet derivative $df(\si)$ in local coordinates,
		and $g_\si^*:X\to Y$ as its adjoint.
		Now we assume there exist $$0<c\le C<\infty,\quad 0\le\al<1$$ s.th. the following holds:
			\begin{align}
	\label{5.4'd}
	&\text{		$df(\s):T_\si\Si\to X$ is injective at every $\sigma\in \Sigma$,}\tag{G1}\\
	\label{5.4d}
	&c\eps^{-\al}\norm{\xi}_Y\le \norm{g_\s\xi}_{Y\to X}\le C\eps^{-\al}\norm{\xi}_Y\quad(\xi\in Y)\tag{G2}.
\end{align}

		Condition \eqref{5.4'} says that the parametrization $f$ in \defnref{M}
		is an immersion, and therefore the manifold $M=f(\Si)$ 
		is non-degenerate.
		Condition \eqref{5.4} is to emphasize that this $\norm{g_\s}_{Y\to X}$ is allowed to
		be large, 
		though not as large as $\eps^{-1}$. 
		It is important  to allow $g_\si$ to be large,
		for the following reason:
		In applications to interface dynamics, 
		$g_\si$ is often a multiplication operator by the gradient of a function that
		has  steep transition layers (e.g. \eqref{5.3}, \eqref{5.2}).
		As such, \eqref{5.4} arises naturally with $\al\ge0$.
		See some concrete examples in \secref{sec:5}.

		In the remaining of this subsection,
		we explain the implication of \eqref{5.4'}-\eqref{5.4} in relation to the geometric structure induced by the 
		evolution \eqref{1}.
		
		Define a  bilinear form
		\begin{equation}
			\label{omega}
			\omega:(u,v)\mapsto \inn{-J^{-1}u}{v}_X,
		\end{equation}
		where either $J=-1$ or $J$ is a symplectic operator satisfying 
		\begin{equation}\label{J}
			J^{-1}=J^*=-J.
		\end{equation}
		If \eqref{1} is a
		gradient flow, then $\omega$ is just the inner product on  $X$. This turns $M$ into a Riemannian manifold.
		If \eqref{1} is Hamiltonian, then $\omega$ is a 
		non-degenerate symplectic form on $X$ and therefore turns $M$ into a symplectic manifold.
		
		Recall that the tangent space $T_\si\Si$ to the Riemannian manifold $\Si$
		can be trivialized 
		as a Hilbert space $Y$.
		Define
		\begin{equation}
			\label{2.1}
			\fullfunction{\J_\s}{Y}{Y}{\xi}{g_\s^*J^{-1}g_\s\xi}.
		\end{equation}
		If $J$ is a symplectic operator, then $\J_\s$ 
		induces a symplectic form on the  tangent bundle $T\Sigma$,
		since
		$$\inn{\J_\s\xi}{\xi}=\inn{g_\s^*J^{-1}g_\s\xi}{\xi}=\inn{J^{-1}g_\s\xi}{g_\s\xi}=0\quad (\xi\in Y).$$
		Moreover,  $\J_\s$ is invertible precisely because of the assumption
		\eqref{5.4'} that  $g_\s$ is injective (or, equivalently, that 
		$f$ is an immersion).
		
		Hence, in both gradient and Hamiltonian cases,
		the operator $\J_\s$ induces a  non-degenerate bilinear  form on $\Sigma$
		through
		$$T_\s\Sigma\times T_\s\Sigma\ni(\xi,\eta)\mapsto \inn{\J_\sigma\xi}{\eta}_Y.$$
		This map is the pull back of \eqref{omega} by the 
		parametrization $f$ (up to a sign). 
		The non-degeneracy of $\J_\s$ is essential in the subsequent sections,
		and this is the most important implication of \eqref{5.4'}.
		This importance was already noted in \cite{MR2189216}.

		Next, we summarize the key implications of \eqref{5.4} into the following lemma:
		\begin{lemma}[estimates on $\J_\si$]
			There hold the following estimates:
		\begin{align}
			\label{5.5}
			\norm{\J_\s}_{Y\to Y}&\ls\eps^{-2\alpha},\\
			\label{5.6}
			\norm{\J_\s^{-1}}_{Y\to Y}&\ls\eps^{2\alpha},\\
			\label{5.8}
			\norm{\E'(\s)}_Y&\ls\eps^{1-\alpha}.
		\end{align}
		Here $\E:\Sigma\to \Rb$ is the pull-back of $E$ by $f$.
		\end{lemma}
		
		\begin{proof}
			Estimate \eqref{5.5} is straight forward. 
		For \eqref{5.8}, we use 
		the identity 
		$$g_\s^*E'(f(\s))=\E'(\s).$$
		This follows from the chain rule.
		
		For \eqref{5.6}, we note that the assumption \eqref{5.4} means the bounded 
		self-adjoint
		operator $\J_\si^*\J_\si:Y\to Y$ 
		has a spectral gap at $0$ of  size $O(\eps^{-4\al})$.
		To see this, we compute
		\begin{align*}
			\inn{\J_\si^*\J_\si\xi}{\xi}&=\inn{(g_\si^* Jg_\si)(g_\si^* J^{-1}g_\si)\xi }{\xi}\\
			&=\inn{(g_\si g_\si^*) J^{-1}g_\si\xi }{J^{-1}g_\si\xi}\\
			&\ge c_1\eps^{-2\al}\norm{J^{-1}g_\si\xi }_Y^2\\
			&\ge c_2\eps^{-4\al}\norm{\xi }_Y^2.
		\end{align*}
		Here $c_1,\,c_2>0$ depend only on the constant $c$ in \eqref{5.4}.
		In the last two inequalities, we use the lower bound from \eqref{5.4} twice,
		together with the fact that $\norm{J(\cdot)}_X=\norm{J^{-1}(\cdot)}_X=\norm{\cdot}_X$
		if $J=-1$ or $J$ is  a symplectic operator  satisfying \eqref{J}.
		By the spectral theorem, we conclude \eqref{5.6} from the lower bound above.
		\end{proof}
		
		\begin{remark}\label{remSi}
			Lastly, we remark that in \defnref{M},
			the closedness assumption on  $\Si$
			naturally arises when $\Si$ represents 
			finitely many small perturbation parameters.
			However, this assumption
			is not essential. For all results in this paper, 
			one can replace this assumption
			by 
			$$\norm{L_\si}_{X\to X}\ls 1,\quad \sup_{\si\in\Si}\left(\norm{f(\si)}_X+\norm{d_\si g_\si}_{Y\to L(Y,X)}\right)\ls 1.$$
			In applications where $\Si$ is unbounded  (e.g. 
			when the moduli represent points in $\Rb^d$), or when $\Si$ is not compact (e.g. when the moduli represent local gauge symmetries), such uniform $O(1)$ estimates as above
			can often be derived
			\textit{ad hoc}.
			See e.g. \cites{MR2094474,MR2097576,MR2232367} for the first case, and \cite{MR2189216} for both cases.
		\end{remark}
		
		\section{Two Key Lemmas}\label{sec:3}
		In this section we 
		prove two key lemmas for the adiabatic theory developed in Sections \ref{sec:G}-\ref{sec:H}.

		Recall that the tangent space 
		$T_{f(\s)}M$ can be trivialized as a subspace of $X$ (see \defnref{M}).
		Define the linear projection $Q_\s$ onto $T_{f(\s)}M$  by 
		\begin{equation}
			\label{2.4}
			\fullfunction{Q_\s}{X}{T_{f(\s)}M\subset X}{\phi}{g_\s \J_\s^{-1}g_\s^*J^{-1}\phi}.
		\end{equation}
		Then $Q_\si^2=Q_\si$ and 
		$$\ran Q_\si=T_{f(\si)}M,\quad \ker Q_\si=(JT_{f(\si)} M)^\perp.$$
		Compare this to the requirements in \eqref{C2}-\eqref{C3}.
		
		Either $J=-1$, or $J$ is a symplectic operator satisfying 
		\eqref{J},
		the operator $Q_\si$ satisfies 
		the identity 
		\begin{equation}\label{2.4'}
			Q_\s^* =J^*Q_\s J.
		\end{equation}
		If $J=-1$, then this implies the projection $Q_\si$ is orthogonal,
		i.e. $Q_\si^*=Q_\si$. If $J$ is symplectic,
		then $Q$ is skew orthogonal with $Q_\si^*=-JQ_\si J$.
		
		Recall
		$$
		\omega:(u,v)\mapsto \inn{-J^{-1}u}{v}_X\quad (u,v\in X)
		$$
		is the bilinear map defined in \eqref{omega}.
		Recall also that if $V\subset X$ is a non-empty subset, and $u$ is an element in $X$,
		then 
		\begin{equation}\label{dist}
			\dist(u,V)=\inf\Set{\norm{u-v}_X:v\in V}.
		\end{equation}
		
		With these definitions at hand, we first  construct a nonlinear 
		projection from
		a tubular neighbourhood around the manifold  $M$ of 
		approximate solitons into the moduli space $\Sigma$. 
		If $u\in X$ lies in this tubular neighbourhood, then
		we call the projection $\s\in \Si$ of $u$ the \textit{moduli} of $u$.
		Indeed, as far as the effective dynamics is concerned, one 
		can think of the vector  $\s$ as a (possibly infinite) tuple of  modulation parameters.
		
		\begin{lemma}[existence of moduli]
			\label{lem2.1}
			Fix two constants $\g, C>0$. There exists $0<\eps_0\ll1$ depending on $\g,\,C$ only s.th.
			the following holds:
			\begin{enumerate}
				\item (Existence of projection) 	There exists an
				open
				neighbourhood $M'\subset X$ around $M$,
				together with a $C^1$ map
				$$S:M'\to \Si,$$
				s.th. for every
				$u\in M'$ and $\si=S(u)$,
				there holds
				\begin{equation}
					\omega(u-f(\s),\phi)=0 \quad (\phi\in T_{f(\s)}M). \label{2.5}
				\end{equation}
				\item (Definite size of $M'$)
				For every $0<\eps\le \eps_0$ and
				$\delta=C\eps^\g$, the set $M'$ can be taken to be
				the tubular neighbourhood
				$$M'=\Set{u\in X:\dist(u,M)<\delta}.$$
				\item (A priori estimate)
				If $u\in M'$ 
				and $u=f(\si_1)+w_1$ for some $\si_1\in \Si$,
				then there holds
				\begin{equation}\label{2.5'}
					\norm{u-S(u)}_X\ls \norm{w_1}.
				\end{equation}
				
			\end{enumerate}
		\end{lemma}
		\begin{remark}
			\label{remS}
			The map $S$ is a projection in the sense that 
			$S(f(\si))=\si$ for every $\si\in \Si$, 
			and $S(f(S(u)))=S(u)$ for every $u\in M'$.
			
			For fixed $\g>0$, the threshold $\eps_0\to 0$ as 
			$C\to \infty$. Hence it is not possible to extend the size of $M'$ indefinitely.
			For fixed $C$, the threshold $\eps_0\to0$ as $\g\to0$ from above.

			Notice that for this lemma, we only need the geometric assumptions from \secref{sec:3.2}. 
		\end{remark}
		\begin{proof}[Proof of the lemma]
			1. 
			First, note that 
			$Q_\s\phi=\phi$ if  $\phi\in T_{f(\s)}M$.
			(One can check this by
			writing $\phi=g_\s\xi$ for some $\xi\in Y$.)
			By this fact, 
			together with the identity \eqref{2.4'}, which holds for both
			gradient and Hamiltonian case,
			we find 
			\begin{align*}
				\om(u-f(\s),\phi)&=\inn{-J^{-1}(u-f(\si))}{Q_\si\phi}\\
				&=-\inn{Q^*J^{-1}(u-f(\si))}{\phi}\\
				&=-\inn{J^*Q_\si(u-f(\si))}{\phi}\quad (\phi\in T_{f(\si)}M).
			\end{align*}
			Hence,  condition \eqref{2.5}
			is satisfied if $Q_\s(u-f(\s))=0$, and	
			for the existence part,
			it suffices to find a map $S:u\mapsto \si$ s.th.
			$$Q_\s(u-S(u))=0.$$
			Consider the map 
			$$\fullfunction{F}{X\times \Sigma}{Y}{(u,\s)}{g_\s^*J^{-1}(u-f(\s))}.$$
			It is clear that if $F(u,\s)=0$, then $Q_\s(u-f(\s))=0$.
			Moreover, 
			if the parametrization map $f$ is $C^2$,
			then $F$ is $C^1$.
			Thus, we proceed to solve the equation
			\begin{equation}\label{F}
				F(u,\s)=0
			\end{equation}
			by Implicit Function Theorem.
			
			Fix a point $\si\in \Si$. 
			The equation \eqref{F} has the trivial solution
			$(f(\s),\s)$. 
			The partial Fr\'echet derivative ${\di_\s F}\vert_{(f(\s),\s)}$ equals to
			$-\J_\s,$ which is 
			invertible as we discussed in \secref{sec:3.2}. 
			Hence, we conclude from Implicit Function Theorem that 
			there exists $\delta=\delta(\s,\eps)>0$
			and a $C^1$ map $S_\s:B_\delta(f(\s))\to \Sigma$ s.th. 
			$F(u,S_\si(u))=0$ for $u\in B_\delta(f(\s))$.
			
			Since $\si$ is arbitrary in the above construction,
			we can patch together all these $S_\si$
			to get an open set $M'\subset X$ containing $M$, together with a $C^1$
			map $S:M'\to X$, s.th. $(u,S(u))$ solves  \eqref{F} for every $u\in M'$.
			
			2. 
			At this point the open set $M'=\bigcup_{\si\in \Si} B_{\delta(\si,\eps)}(f(\si))$.
			Now we claim in fact $\delta$ can be made independent of $\s$. This is essential for our purpose, because
			we would like $M'$ to contain
			a definite volume  for a  flow to fluctuate.
			
			Fix a point $\si\in\Si$. Write 
			$$A_w:={d_\s F}\vert_{(w+f(\s),\s)},\quad V_w:=A_w-A_0.$$
			Then $$A_0=-\J_\s,\quad V_w={(d_\s g_\s^*)\del{\cdot}}\vert_{J^{-1}w}.$$
			We recall from the proof of Implicit Function Theorem  
			(e.g. \cite{MR1336591}*{Sec. 2}) 
			that the size of $\delta>0$ from
			the above construction is determined by the following condition:
			For every $w\in B_\delta(f(\s))$, there hold
			\begin{align}
				&\label{2.4.00}
				A_w\text{ is invertible},\\
				&\label{2.4.0}\norm{A_w^{-1}}_{Y\to Y}\le2\norm{A_0^{-1}}_{Y\to Y},\\
				&\label{2.4.1}
				\norm{F(w+f(\s),\s)}_{Y}\le \frac{\delta_0}{4\norm{A_0^{-1}}_{Y\to Y}},
			\end{align}
			where $\delta_0>0$ is chosen so that  for the remainder
			$$R(w,\xi):=F(w+f(\s),\s+\xi)-F(w+f(\s),\s)-\di_\s F(w+f(\s),\s)\xi,$$
			the following conditions hold
			for every $\xi\in B_{\delta_0}(\s)$ and $u\in B_\delta(f(\s))$:
			\begin{align}
				&\norm{R(w,\xi)}_{Y}\le  \frac{\delta_0}{4\norm{A_0^{-1}}_{Y\to Y}},\label{2.4.2}\\ 
				&\norm{\J_{\s+\xi}-\J_\s}_{Y\to Y}\le\frac{1}{4\norm{A_0^{-1}}_{Y\to Y}}.\label{2.4.3}
			\end{align}
			Note that the r.h.s. of \eqref{2.4.0}-\eqref{2.4.3} are independent of $\s$ 
			by the uniform estimates for $\J_\s^{-1}=-A_0^{-1}$ from \eqref{5.6}.
			
			Fix any $c_0>0$,
			and recall $\g>0$ is given.
			We claim conditions \eqref{2.4.2}-\eqref{2.4.3} are satisfied for  
			$\delta_0=c_0\eps^\g$ and all sufficiently small $\eps$.  
			
			To get \eqref{2.4.2}, one uses the fact that $\norm{R(w,\xi)}_{Y}=o(\norm{\xi}_{Y})=o(\delta_0)$, since $R$
			is the super-linear remainder of the expansion of the $C^1$ map $F$ in $\s$.
			At this point we need a constraint  $\eps\le \eps_0$ for 
			some $\eps_0>0$ depending on $\g$.
			This constraint $\eps_0\to 0$ as $\g\to 0$ from above.
			
			To get \eqref{2.4.3}, one uses the continuity of the map
			$\s \mapsto \J_\s \in L(Y,Y)$ and the estimate \eqref{5.6},
			which imply that \eqref{2.4.3} holds so long as $\delta_0=o(1)$.
			At this point we need another constraint $\eps_0=o(c_0^{-1})$ as $c_0\to\infty$.
			Importantly, 
			this constraint implies the size of $\delta$ 
			cannot be made indefinite.
			
			Next, we claim \eqref{2.4.00}-\eqref{2.4.1} are satisfied for 
			$\delta=c_1\delta_0$ with some fixed $c_1$ independent of $c_0$
			and all $\eps\le \eps_0(c_0)$. 
			This, together with the
			choice of $\delta_0$ above, confirms the 
			claim about the size of $\delta$.
			
			Indeed, with the choice $\delta_0=c_0\eps^\g$, condition
			\eqref{2.4.1} is satisfied  if and only if
			\begin{equation}\label{2.4.1'}
				\norm{A_0^{-1}}_{Y\to Y}\norm{F(w+f(\s),\s)}_{Y}\le \frac{c_0\eps^\g}{4}.
			\end{equation}
			By the uniform estimate for 
			$g_\s^*$ and $\J_\s^{-1}$, we find that l.h.s. of this expression
			can be bounded from above by $c_2\eps^\al\norm{w}_X$ for some $c_2>0$ depending
			only on the constants $c,\,C$ in \eqref{5.4}.
			Thus, \eqref{2.4.1'} holds with the choice
			$\delta=c_0\eps^\g/(4c_2)^{-1}.$

			Next, by elementary perturbation theory, since 
			$A_0$ is invertible, it follows that condition \eqref{2.4.00} is satisfied 
			so long as 
			\begin{equation}\label{V}
				\norm{V_w}_{Y\to Y}\le\frac{1}{2}\norm{A_0^{-1}}_{Y\to Y}^{-1}=\frac{1}{2}\norm{\J_\s^{-1}}^{-1}_{Y\to Y}.
			\end{equation}
			By \eqref{5.6}, we have $\norm{\J_\s^{-1}}^{-1}_{Y\to Y}\ge c_3$
			for some $c_3>0$ depending on the implicit constant in \eqref{5.6} only.
			By the condition that $f$ is $C^2$
			and $\Si$ is closed (see \defnref{M} as well as \remref{remSi}),
			we have a uniform $O(1)$ bound on the linear map $\xi\mapsto d_\si g_\si^*(\xi)\in L(X,Y)$.
			By this,  together with the definition of $V_w$, we 
			conclude $\norm{V_w}_{Y\to Y}\le c_4\norm{J^{-1}w}_X=c_4\norm{w}_X$
			for some $c_4>0$ depending on $f$ only.
			Thus,  \eqref{V} is satisfied 
			if $\delta= c_3/(2c_4).$

			Lastly, referring to the Neumann series for the inverse 
			$$A_w^{-1}=\sum_{n=0}^\infty A_0^{-1}\del{-V_w A_0^{-1}}^n 	,$$ 
			we find 
			$$\norm{A_w^{-1}}_{Y\to Y}\le \frac{\norm{A_0^{-1}}_{Y\to Y}}{1-\norm{V_w}_{Y\to Y}\norm{A_0^{-1}}_{Y\to Y}}.$$
			With the previous choice $\delta= c_3/(2c_4)$, we conclude \eqref{2.4.0} from this and \eqref{V}.

			This proves the claim about the size of $M'$, with the choice 
			$$\delta= \eps^\g \min \Set{\frac{c_3}{2c_4},\,\frac{c_0}{4c_2}},$$
			which is valid for arbitrary fixed $\g,\,c_0>0$ and all $0<\eps\le \eps_0(\g,c_0)$.
			
			3.
			Lastly, we establish the estimate \eqref{2.5'}.
			
			Suppose $u\in M'$ 
			and $u=f(\si_1)+w_1$ for some $\si_1\in \Si$.
			Let $S$ be the nonlinear projection constructed above,
			and let $w:=u-f(S(u))$.
			Then 
			\begin{equation}\label{2.5.1}
				w=f(\si_1)+w_1-f(S(u)).
			\end{equation}
			Estimate \eqref{2.5'} is equivalent to 
			the bound $\norm{w}_X\ls \norm{w_1}_X$,
			which we prove below.
			
			Consider the expansion
			\begin{equation}\label{2.5.2}
				\begin{aligned}
					f(S(u))=&f(S(f(\si_1)+w_1))\\
					=&f(\si_1+dS(f(\si_1))w_1+o_{\norm{\cdot}_X}(w_1))\\
					=&f(\si_1)+g_{\si_1}(dS(f(\si_1))w_1+o_{\norm{\cdot}_X}(w_1))\\
					&+o_{\norm{\cdot}_Y}(dS(f(\si_1))w_1).
				\end{aligned}
			\end{equation}
			The second line is valid since $S$ is $C^1$ and $S(f(\si_1))=\si_1$.
			The third line 
			is valid since $f$ is $C^2$.
			In view of \eqref{2.5.1} and \eqref{5.4}, 
			it remains to find a uniform $O(\eps^\al)$ estimate on the linear operator
			$dS(u):X\to Y$.

			Differentiating the equation $F(u,S(u))=0$, we find 
			$$0=d_uF(u,S(u))=\di_u F(u,S(u))+\di_\si F(u,S(u))dS(u).$$
			This implies $dS(u)=-(\di_\si F(u,S(u)))^{-1}\di_u F(u,S(u))=-\J_{S(u)}^{-1} g_{S(u)}^*J^{-1}.$ By \eqref{5.4} and \eqref{5.6}, we conclude 
			$\norm{dS(u)}_{X\to Y}\ls \eps^\al$. Plugging this into \eqref{2.5.2},
			and using \eqref{5.4},
			we find $f(S(u))=f(\si_1)+O_{\norm{\cdot}_X}(w_1)$.
			Hence, the desired estimate \eqref{2.5'} follows from \eqref{2.5.1}.
			
			This completes the proof.
		\end{proof}

		By \lemref{lem2.1}, 
		if $u_t$
		is a path in $X$ with 
		$\dist(u_t,M)\ls \eps^\g$ for some $\g>0$
		and sufficiently small $\eps>0$, then there holds 
		the unique decomposition
		\begin{equation}\label{2.5.3}
			u_t=f(\si_t)+w_t\quad \text{s.th.} \quad \si_t\in \Si,\; Q_{\si_t}w=0.
		\end{equation}

		Moreover, the choice $\s_t$ is optimal in the sense that
		the (skew) orthogonality condition \eqref{2.5} is satisfied.
		To see \eqref{2.5} is a natural condition for optimality,
		we note that if $J=-1$, then \eqref{2.5} means $w\perp T_{f(\si)}M$
		and this guarantees $f(\s_t)$ is the closest 
		path in $M$ to $u_t$. 
		In the Hamiltonian case,
		in the presence of continuous symmetry,
		the 
		skew orthogonality condition is also customarily
		used 
		to derive the modulation equations
		for solitary wave dynamics, see e.g. \cites{MR1071238,MR1170476,MR2094474,MR2232367,MR2351368}.
		
		\begin{remark}\label{remD}
			In the sequel we will use the a priori estimate \eqref{2.5'} as follows:
			Suppose $u\in X$ satisfies  $\dist(u_,M)\le C \eps^\g$
			for some $\g,\,C>0$ and sufficiently small $\eps$.
			Then 
			there is $\si_*\in \Si$ s.th.
			$\norm{u-f(\s_*)}_X\le 2C\eps^\g$ by definition \eqref{dist}.
			Now, if $\si=S(u)$ is the moduli associated to $u$
			and $w:=u-f(\si)$, then
			applying 
			the a priori estimate \eqref{2.5'}
			yields $\norm{w}_X\ls\norm{u-f(\s_*)}_X\le C'\eps^\g$
			for some $C'\ge2C$.
			This fact allows us to keep track only the leading
			order term in $w$ in the derivation of various remainder estimates below,
			knowing only that $u_t$ is close to some point in $M$.
		\end{remark}
		
		Recall $\E:\Si\to \Rb$
		is the pull-back of the energy functional $E$ by  $f$.
		If  $\s$ is the moduli of $u$, then $\E(\s)$ is the effective 
		energy of the latter. Hence,  if $u_t$ can be decomposed as \eqref{2.5.3}, then heuristically, 
		one expects  the effective dynamics governing 
		the motion of the moduli $\si=\s_t$ to be
		\begin{equation}
			\label{eff}
			\di_t\s= \J_\s^{-1}\E'(\s) .
		\end{equation}
		In particular, the energy property of \eqref{eff} (i.e. dissipative or
		conservative) agrees with that of \eqref{1}.

		In the next lemma, we justify the heuristic  choice of \eqref{eff} 
		as the effective dynamics for a full flow 
		$u_t$ solving \eqref{1}, \textit{assuming} $u_t$ stays
		uniformly close to $M$.
		We drop this assumption in the next two sections,
		and we  will show it suffices to have $u_t$ near $M$ only at $t=0$
		for the approximation \eqref{eff} to be  valid 
		globally in the gradient flow case, and 
		up to some large time in the Hamiltonian case.
		
		\begin{lemma}
			\label{lem2.2}
			Let $0<T\le\infty$.
			Let $u_t,\,0\le t<T$ be a solution to \eqref{1}.
			Suppose $\dist(u_t,M)\ls\eps^\g$ for all $t\le T$ and some $\g>0$.
			Write $u_t=f(\si_t)+w_t$ as in \eqref{2.5.3}. Then there holds the following uniform estimate for all $t\le T$:
			\begin{equation}
				\label{2.6}
				\norm{\di_t\s-\J_\s^{-1}\E'(\s)}_{Y}\ls\eps^{1+\alpha}\norm{w}_X.
			\end{equation}
		\end{lemma}
		\begin{proof}
			1. Expand \eqref{1} as
			\begin{equation}
				\label{3.5}
				\di_tv+\di_tw=J(E'(v)+L_\s w+N_\s(w)),
			\end{equation}
			where $L_\s$
			is the linearized operator at $v_t:=f(\s_t)$, and $N_\s(w)$
			defined by this equation. 
			This expansion holds by the $C^2$ regularit of $E$.
			
			Recall $Q_\s$ is the projection onto $T_vM$ defined in \eqref{2.4}.
			Applying $Q_\s$ to both sides of \eqref{3.5}, we have
			\begin{equation}
				\label{3.7}
				\di_tv-Q_\s JE'(v)=Q_\s(JL_\s w -\di_tw+JN_\s(w)).
			\end{equation}
			
			Consider the identity
			$$\J_\s^{-1}g_\s^*J^{-1}(\di_tv-Q_\s JE'(v))=\di_t\s-\J_\s^{-1}\E'(\s).$$
			To verify this, one uses two facts that follow readily from the chain rule:
			$$\di_tv=g_\s\di_t\s,\quad g_\s^*E'(f(\s))=\E'(\s).$$
			Thus by the uniform estimates for $g_\s^*$ and $\J_\s^{-1}$, we have
			\begin{equation}
				\label{3.8}
				\norm{\di_t\s-\J_\s^{-1}\E'(\s)}_{Y}\ls\eps^{\alpha}\norm{\di_t v-Q_\s JE'(v)}_X.
			\end{equation}
			
			2. Consider now the r.h.s. of \eqref{3.7}. 
			These three terms can be bounded respectively as follows:
			\begin{align}
				\norm{Q_\s JL_\s w}_{X}&\ls\eps\norm{w}_{X}\label{3.9},\\
				\norm{Q_\s\di_tw}_{X}&\ls\eps^{-\al}\norm{\di_t\s}_{Y}\norm{w}_{X},\label{3.10}\\
				\norm{Q_\s JN_\s(w)}_{X}&\ls\norm{w}_{X}^2.\label{3.11}
			\end{align}
			In all these three inequalities we use the uniform bound $\norm{Q_\s}_{X\to X}\ls 1$.
			
			For \eqref{3.9} we need  the identity 
			\begin{equation}\label{id}
				\abs{\inn{Q_\s JL_\s w}{w'}}= \abs{\inn{w}{L_\s Q_\s Jw'}}.
			\end{equation}
			In both gradient and Hamiltonian case, we have $Q_\si J=JQ^*$
			by  \eqref{2.4'}, and \eqref{id} follows from here.
			
			By \eqref{id}, we find
			\begin{equation}\label{3.200}	
				\begin{aligned}
					\abs{\inn{Q_\s JL_\s w}{w'}}&= \abs{\inn{w}{L_\s Q_\s Jw'}}\\&\le \norm{L_\s Q_\s }_{X\to X}\norm{w}_{X^1}\norm{w'}_{X}
					\quad (w,w'\in X).
				\end{aligned}
			\end{equation}
			Plugging $w'=Q_\s J_\s Lw$ into \eqref{3.200}, we get 
			$$\norm{Q_\s JL_\s w}_{X}\le \norm{L_\s Q_\s }_{X\to X}\norm{w}_{X}\ls\eps\norm{w}_{X}.$$
			The last inequality follows from the approximate zero mode property
			\eqref{C3}.

			Next, for \eqref{3.10}, we use
			the construction from \lemref{lem2.1},
			which ensures the remainder $w=u-v$ satisfies 
			$Q_\si w=0$. We note that precisely at this point
			we use this optimal construction  in an essential way. See a discussion in \remref{remW} below. 
			
			Indeed, differentiating 
			$Q_\s w=0$ w.r.t. $t$, we find
			\begin{equation}\label{3.10.1}
				0=\di_t(Q_\s w)=(\di_t Q_\s)w+Q_\s\di_tw=(d_\s Q_\s\di_t\s)w+Q_\s\di_tw.
			\end{equation}
			Here $d_\s Q_\s$ is an operator from $Y$ to the space of linear
			operators $L(X, X)$.  
			Geometrically,
			since $Q_\s$ is the projection onto the tangent space $T_{f(\si)} M$,
			the operator $d_\si Q_\si$ is the Weingarten map (or shape operator),
			and therefore the bound on $d_\si Q_\si$ depends only on
			the curvature on $M$.
			Since the map $f:\Si\to M\subset X$ is a $C^2$ immersion of
			a closed manifold $\Si$ (see also \remref{remSi}),
			we find the uniform estimate
			$$\norm{d_\s Q_\s}_{Y\to L(X, X)}\ls 
			\sup_{\si\in\Si} \del{\norm{g_\si}_{Y\to X}+\norm{d_\si g_\si}_{Y\to L(Y,X)}}\ls
			\eps^{-\al}.$$ 
			Plugging this into \eqref{3.10.1} gives \eqref{3.10}.

			Lastly, 	\eqref{3.11} follows from the remainder estimate $N_\s(w)=o(\norm{w}_X^2)$, since $E$ is $C^2$.
			
			3. Plugging \eqref{3.9}-\eqref{3.11} to \eqref{3.7}-\eqref{3.8} gives 
			\begin{equation}
				\label{2.7}
				\norm{\di_t\s-\J_\s^{-1}d_\s E(f(\s))}_{Y}\ls \del{\norm{\di_t\s}_Y
					+\eps^{1+\alpha} }\norm{w}_X.
			\end{equation}
			Here note that
			as we discussed in \remref{remD}, we can absorb higher order terms in $w$ 
			into the first order ones.

			Now
			we want to estimate $\norm{\di_t\s}_Y$ in the r.h.s. of \eqref{2.7} 
			at the order of $O(\eps^{1+\al})$,
			whence the claim  \eqref{2.6} follows.
			
			Applying the reverse triangle inequality to the l.h.s. of \eqref{2.7}, we find two cases.
			If $\norm{\di_t\s}_Y< \norm{\J_\s^{-1}\E'(\s)}_Y$,
			then $\norm{\di_t\s}_Y\ls\eps ^{1+\alpha}$ by \eqref{5.6}-\eqref{5.8}. Otherwise, if
			$\norm{\di_t\s}_Y\ge \norm{\J_\s^{-1}\E'(\s)}_Y,$
			then \eqref{2.7} implies
			\begin{equation}\label{2.7'}
				\norm{\di_t\si}_Y\le \norm{\J_\s^{-1}\E'(\s)}_Y + C\del{\norm{\di_t\s}_Y\norm{w}_X
					+\eps^{1+\alpha} \norm{w}_X},
			\end{equation}
			where $C>0$ is independent of $\eps$ and time.
			So long as  \begin{equation}\label{2.7.1}
				\norm{w}_X= o(1),\quad 0<\eps\ll1,
			\end{equation}
			we can transpose the second term in the r.h.s. of \eqref{2.7'} to obtain
			\begin{equation}\label{2.7.2}
				\frac{1}{2}\norm{\di_t\si}_Y\le \norm{\J_\s^{-1}\E'(\s)}_Y + C
				\eps^{1+\alpha} \norm{w}_X.
			\end{equation}
			The Ansatz \eqref{2.7.1} holds since by assumption $\dist(u_t,M)\ls \eps^\g$
			for some 	 $\g>0$ and $\eps\ll1.$
			This implies 
			$\norm{w}_X=O(\eps^\g)=o(1)$ as we explained in \remref{remD}.
			
			From \eqref{2.7.2}  we conclude
			\begin{equation}\label{2.8}
				\norm{\di_t \s}_Y\ls\eps^{1+\alpha}+\eps^{1+\alpha}\norm{w}_X\ls \eps^{1+\al}.
			\end{equation}
			Hence, in both cases we have shown the r.h.s. of \eqref{2.7} is of the order $O(\eps^{1+\al}\norm{w}_X)$. Thus  \eqref{2.6} is proved.
		\end{proof}
		
		\begin{remark}\label{remW}
			Here we would like to remark on the estimate \eqref{3.10}.
			Indeed, it is not in general possible to estimate
			the the full velocity $\di_tw$, because this 
			fluctuation field, however small, may vary rapidly, especially in the Hamiltonian space
			due to acceleration.
			However, an estimate on the projection $Q_\si \di_tw$ is possible because of the identity \eqref{3.10.1}.
			This identity is an important  consequence of the (skew)
			orthogonality condition \eqref{2.5}, as 
			\eqref{3.10.1} relates the 
			tangential velocity $Q_\si \di_tw$ 
			to the velocity of moduli, $\di_t\si$.
			The latter is small up to a large time, so long as initially $\dist(u_0,M)\ll\eps$, as we show in the next sections.

		\end{remark}

		\section{Effective Dynamics for Gradient Flow}\label{sec:G}
		In this section we consider \eqref{1} with $J=-1$	  on 
		the tangent bundle $TU$. In this case the evolution reads
		\begin{equation}
			\label{3.1}
			\di_tu=-E'(u).
		\end{equation}
		This turns \eqref{1} into 
		the gradient flow of $E$. We show any flow starting near 
		the manifold of approximate soliton $M$ can be
		approximated uniformly for all time by a gradient flow
		of the effective energy $\E$ on the moduli space $\Sigma$. 
		Then we derive a converse of this as a corollary.
		
		\begin{theorem}
			\label{thm1}
			Fix any $0<\eps\ll1$. There exists $0<\delta\ll\eps$
			s.th.
			the following holds:
			Let $M$ be the manifold of approximate solitons
			as in \defnref{M}.
			Let $u_0\in X$ be an initial configuration s.th.
			$\dist(u_0,M)\le\delta$. Let $u_t$ be the flow
			generated by $u_0$ under \eqref{3.1}.
			\begin{enumerate}
				\item (A priori estimate) For all $t\ge0$, there holds  
			\begin{equation}
				\label{3.2}
				\dist(u_t,M)\ls \eps.
			\end{equation}
			
			\item (Effective dynamics) Moreover, 
			the decomposition \eqref{2.5.3} for $u_t$ is valid for all time, and the moduli
			$\s\equiv \si_t:=S(u_t)$ satisfies the following effective dynamics:
			\begin{equation}
				\label{3.3}
				\di_t\s=-(g_\s^*g_\s)^{-1}\E'(\s)+O_{\norm{\cdot}_Y}(\eps^{2+\al}).
			\end{equation}
			\end{enumerate}
		\end{theorem}
		\begin{remark}
			
			The remainder in \eqref{3.3} is of lower order by
			\eqref{5.8}.
		\end{remark}
		\begin{proof}[Proof of Theorem \ref{thm1}]

			1.  To begin with, note that by the continuity of the flow \eqref{3.1},
			if $\delta\ll\eps$, then 
			there exists some (possibly small) $0<T_1\le \infty$ s.th.
			\eqref{3.2} holds for $t<T_1$.
			This gives the decomposition 
			\begin{equation}\label{decomp}
				u_t=v_t+w_t,\quad v_t:=f(\si_t),\quad \si_t:=S(u_t)
			\end{equation}
			as in \eqref{2.5.3}, which is valid for $0\le t<T_1$.
			Here $w_t$ is defined by the relation \eqref{decomp},
			i.e. $w_t:=u_t-v_t=u_t-f(S(\si_t))$. 
			
			The claim now is that we have the a priori estimate 
			\begin{equation}\label{A1}
				\norm{w_t}_X\le C(1+e^{-\g t})\eps\quad(t\le T_1)
			\end{equation}
			for some $\g,\,C>0$ independent of $t$ and $T_1$.
			If this holds,  then since the constant 
			$C$ is independent of time, a standard blow-up argument
			yields $T_1=\infty$, and \eqref{3.2} follows
			since by definition \eqref{dist}, we have 
			$\dist(u_t,M)\le \norm{w_t}_X$.
			\lemref{lem2.1} then guarantees the validity
			of the decomposition \eqref{decomp}
			for all time, and 
			 the remainder estimate in \eqref{3.3} follows from \eqref{2.6} and \eqref{A1}.

			 Hence, the theorem
			 is proved once we establish \eqref{A1}.
			 
			To this end,
			we derive a differential inequality for the 
			function
			\begin{equation}\label{Lyap}
				t\mapsto \frac{1}{2}\inn{L_{\si(t)} w(t)}{w(t)},
			\end{equation}
			which accounts for  most of the energy dissipation.
			We will show this quadratic form is approximately a
			Lyapunov functional along \eqref{3.1}. 
			Then by  the
			coercivity condition \eqref{C2},
			this approximately monotone quantity  controls $\norm{w}_X$,
			since by the orthogonality condition \eqref{2.5},
			the fluctuation field $w\in \ker Q_\si=(T_{f(\si)}M)^\perp$.

			2.  We now study the 
			quantity $\frac{1}{2}\inn{L_\s w}{w}.$
			Compute
			\begin{equation}
				\label{16}
				\begin{aligned}
					\frac{1}{2}\od{\inn{L_\s w}{w}}{t}&=\inn{\di_tw}{L_\s w}+\frac{1}{2}\inn{(\di_t L_\s)w}{w}\\&=\inn{-\di_tv-(E'(v)+L_\s w+N_\s(w))}{L_\s w}+\frac{1}{2}\inn{(\di_t L_\s)w}{w}.
				\end{aligned}
			\end{equation}
			Here we have used the expansion \eqref{3.5}.
			We bound the two inner products in the last line of \eqref{16}. 
			
			To bound the second one, we note two things:
			First, there holds the identity $\di_t L_\s=(d_\s L_\s )\di_t\s$ by the chain rule.
			Second, we have a uniform
			bound on $d_\s L_\s:Y\to L(X,X)$ of the order $O(\eps^{1-\al})$.
			To see this, we compute $d_\si L_\si=d_\si E''(f(\si))=d_v\vert_{v=f(\si)}E''(v)df(\si)=L_\si g_\si$.
			For the last equality, we note
			that for a linear map,  the Fr\'echet derivative is itself.
			Since $g_\si$ maps into (the trivilization of) $T_{f(\si)}M$,
			the claimed uniform bound on $d_\si L_\si$ follows from the assumptions \eqref{C3}
			and \eqref{5.4}.

			From the preceding discussion, we conclude the following estimate for the second
			term in the last line of \eqref{16}:
			\begin{equation}
				\label{19}
				\abs{\inn{(\di_t L_\s)w}{w}}\ls\eps^{1-\al}\norm{\di_t\s}_Y\norm{w}_X^2.
			\end{equation}
			We note that the approximate zero-mode property \eqref{C3} is used
			here in a crucial way to derive this estimate.

			Now we claim the following three estimates hold:
			\begin{align}
				\inn{-\di_tv-E'(v)}{L_\s w}&\ls(\eps+\norm{w}_X)\norm{w}_X^2+\eps\norm{w}_X,\label{19.1}\\
				-\inn{N_\s(w)}{L_\s w}&\ls\norm{w}_X^3,\label{19.3}\\
				-\inn{L_\s w}{L_\s w}&\le -\beta'\norm{w}_X^2\text{ for some fixed }\beta'>0.\label{19.2}
			\end{align}
			For all these estimates we need a uniform bound on $\norm{L_\si}_{X\to X}$.
			Recall in \defnref{M} we assume the moduli space $\Si$ to be closed.
			Moreover, we assume
			the map $f$ and the energy functional $E$ are both $C^2$.
			These facts imply
			that the map $\si\mapsto L_\si$ is continuous and  bounded on $\Si$. 
			Consequently, there exists some fixed $C>0$ s.th.
			\begin{equation}\label{20'}
				\norm{L_\s}_{X\to X}\le C.
			\end{equation}
			See also \remref{remSi}.

			For \eqref{19.1}, we recall equation \eqref{3.7} and 
			the estimates \eqref{3.9}-\eqref{3.11}, \eqref{2.8}
			derived in \lemref{lem2.2}.
			Rearranging \eqref{3.7}, we conclude from these estimates and condition \eqref{C1} that there
			holds the velocity bound
			\begin{align*}
				\norm{\di_tv}_X&\ls \eps+(\eps+\eps^{-\al}\norm{\di_t\si}_X+\norm{w}_X)\norm{w}_X\quad \text{by \eqref{C1} and \eqref{3.9}-\eqref{3.11}}\\
				&\ls \eps +(\eps+\norm{w}_X)\norm{w}_X\quad \text{by \eqref{2.8}}.
			\end{align*}
			This and another application of \eqref{C1} gives \eqref{19.1}.

			For \eqref{19.3}, we use 
			the nonlinear estimate $\norm{N_\s(w)}_X\ls\norm{w}_X^2$,
			which follows from the $C^2$ regularity of the energy functional $E$.
			
			For \eqref{19.2}, we recall that the orthogonality 
			condition from \lemref{lem2.1} in the gradient case ensures
			$w\in (T_vM)^\perp$. By this fact, \eqref{19.2}
			follows from  the stability condition \eqref{C2},
			and the  constant $\beta'$
			depends on the gap size $\beta$ in \eqref{C2} only.
			(In general we have $\beta'\le\beta^2$.)
			This spectral gap condition is precisely used here to
			get the bound \eqref{19.2}, which is the central estimate in what follows.
			
			3. Combining \eqref{19.1}-\eqref{19.3}, we get an estimate 
			\begin{equation}
				\label{18}
				\inn{-\di_tv-(E'(v)+L_\s w+N_\s(w))}{L_\s w}\ls(\norm{w}_X+\eps-\beta')\norm{w}_X^2.
			\end{equation}

			Plugging \eqref{19} and \eqref{18} into \eqref{16}, we find
			\begin{equation}
				\label{20}
				\frac{1}{2}\abs{\od{\inn{L_\s w}{w}}{t}}\ls(\norm{w}_X+\eps-\beta')\norm{w}_X^2+\eps\norm{w}_X.
			\end{equation}
			
			Now, by \eqref{20}, so long as $\eps\le \beta/4$,
			we can find $\g>0$ depending on the constant $C$ in \eqref{20'} and 
			$\beta'$ in \eqref{19.2} only s.th. 
			\begin{equation}
				\label{21}
				(\tfrac{d}{dt}+\gamma)\inn{L_\s w}{w}\ls(\norm{w}_X-\beta'/2)\norm{w}_X^2+\eps\norm{w}_X.
			\end{equation}
			
			At this point we make the Ansatz
			\begin{equation}
				\label{A}
				\norm{w}_X\le \beta'/2.
			\end{equation}
			Shrinking $\epsilon$ if necessary, 
			we can ensure that 
			this Ansatz holds at least locally
			for  $t\le T_1$.

			So long as \eqref{A} holds,  we can drop the first term in the r.h.s.
			of \eqref{21}, and then multiply both side by $e^{\gamma t}$ to get 
			\begin{equation}
				\label{21'}
				\od{}{t}\del{e^{\gamma t}\inn{L_\s w}{w}}\ls \eps e^{\gamma t}\norm{w}_X.
			\end{equation}
			
			Integrating \eqref{21'}, and then dividing the integration factor, we find 
			\begin{equation}
				\label{22}
				\begin{aligned}
					\inn{L_\s w}{w}&\ls e^{-\gamma t}\inn{L_{\s_0} w_0}{w_0}+\eps M(t)\\
					&\ls e ^{-\gamma t}\norm{w_0}_X^2+\eps M(t)
					\quad \del{M(t):=\sup_{t'\le t}\norm{w_{t'}}_X}.
				\end{aligned}
			\end{equation}
			Together with the coercivity condition \eqref{C2}, we find
			\begin{equation}
				M(t)\le C_1\del{e^{-\gamma t}\norm{w_0}_X+\eps },\label{25}
			\end{equation}
			where $C_1>0$ depends on the spectral gap from \eqref{C2}
			and is independent of $t$ and $T_1$. 
			
			As we discussed in \remref{remD}, if we now choose $\delta\ll\eps$,
			then $\norm{w_0}_X\le\eps$. This, together with \eqref{25} above, 
			implies the key a priori estimate \eqref{A1}.
			We also conclude from here that Ansatz \eqref{A} holds so long 
			as $\eps$ is sufficiently small.
			This completes the proof.

		\end{proof}

		\begin{corollary}[Converse of \thmref{thm1}]\label{cor1}
			Fix any large $T\gg1$ and $0<\g<1$.
			There exist $c>0$ independent of $T$, and 
			$0<\eps_0\ll1$ depending on $\g,\,T$ only,
			s.th. the following holds:
			Let $0<\eps\le\eps_0$.
			Let $\s_0\in\Sigma$. Let $\s_t\in\Sigma$ be the flow generated by $\s_0$
			under the effective dynamics \eqref{eff}. Then there exists a solution $u_t$ to \eqref{3.1}
			s.th. for all $\epsilon^\g t\le T$, there holds
			\begin{equation}
				\label{23}
				\norm{u_t-f(\sigma_t)}_X\le c \epsilon^{(1-\g)/2}.
			\end{equation}
		\end{corollary}
		\begin{remark}
			Note that  this result holds only on a long finite interval, and the remainder 
			in \eqref{23} tends to $0$ as $\eps\to0$.
			
			For fixed $\g>0$, the threshold $\eps_0\to 0$ as $T\to \infty$.
			For fixed $T$, the threshold $\eps_0\to 0$
			as $\g\to 1$ from below. 
		\end{remark}
		\begin{proof}[Proof of the Corollary]
			Let $u_0=f(\s_0)$. Let  $u_t$ be the flow generated by $u_0$ under
			\eqref{1}. 
			We claim this flow $u_t$ satisfies \eqref{23} for $\eps^\g t\le T$.
			
			Indeed, by \thmref{thm1}, there exists a flow $\tilde{\s}_t\in\Sigma$
			s.th. 
			\begin{align}
				&\norm{u_t-f(\tilde{\s}_t)}_X\le c_1 \epsilon,\label{24.1}\\
				&\norm{\di_t\tilde{\s}-\J_{\tilde{\si}}^{-1}\E'(\tilde{\s})}_Y\le c_2\eps^{2+\al}.\label{24.2}
			\end{align}
			Here $c_1,\,c_2>0$ are some constants independent of $\eps,\,T,$ and $c$.
			
			By the uniqueness of moduli (see \lemref{lem2.1}), at the initial time we have $\sigma_0=\tilde{\s}\vert_{t=0}$.
			Using this, and integrating \eqref{24.2}, we find $$\begin{aligned}
				\norm{\tilde{\s}_t-\s_t}_Y&\le \eps^{-\g}T\sup_{\eps^\g \tau\le T}
				(\norm{\di_t\tilde{\s}_\tau-\J_{\tilde{\si_\tau}}^{-1}\E'(\tilde{\s_\tau})}_Y+\norm{\J_{\tilde{\si}_\tau}^{-1}\E'(\tilde{\s}_\tau)-\J_{\si_\tau}^{-1}\E'(\s_\tau)}_Y)\\
				&\le T(c_2\eps^{2+\al-\g}+c_3\eps^{1+\al-\g}),
			\end{aligned}$$
			so long as $\eps^\g t\le T$. In the last line we have used 
			\eqref{24.2}, and the constant $c_3$ depends on the implicit
			constants in \eqref{5.6} and \eqref{5.8} only.
			
			By the uniform estimate \eqref{5.4}, it follows that 
			$$
			\norm{f(\tilde{\s_t})-f(\s_t)}_X\le C \eps^{-\al}\norm{\tilde{\s}_t-\s_t}_Y\le CT(c_3\eps^{1-\g}+c_2\eps^{2-\g})\quad (\eps^\g t\le T).
			$$
			Here $C>0$ is the constant from \eqref{5.4}.
			
			Applying triangle inequality to \eqref{24.1}, 
			we find
			\begin{equation}\label{24.3}
				\norm{u_t-f(\sigma_t)}_X\le CT(c_3\eps^{1-\g}+c_2\eps^{2-\g}) + c_1\eps \quad (\eps^\g t\le T).
			\end{equation}
			Recall that $T\gg1$ is fixed, the exponent $0<\g<1$, and these $c_1,c_2,c_3,C$
			are all absolute constants.
			Hence, the leading term in \eqref{24.3} is of the order 
			$O(\eps^{1-\g})$, and  from \eqref{24.3} we can choose $\eps_0$ depending on $T$ only
			to conclude \eqref{23} for every $0<\eps\le\eps_0$.
		\end{proof}

		\section{Effective Dynamics for Hamiltonian System}\label{sec:H}
		In this section we consider 
		the Hamiltonian system
		\begin{equation}
			\label{Ham}\di_tu=J E'(u),\quad J^{-1}=J^*=-J.
		\end{equation}
		One essential difference in the analysis of \eqref{Ham} from that of the
		gradient flow \eqref{3.1}
		is the following: For \eqref{Ham}, we do not have a natural
		Lyapunov-type functional that bounds the fluctuation field,
		such as the quadratic form
		\eqref{Lyap}.
		Indeed, the decay property of \eqref{Lyap} is 
		ultimately due to the energy dissipation for gradient flows.
		
		There are two important classes of evolutions of the form \eqref{Ham} that  arise from physics:
		
		1. Let $\Omega\subset \Rb^d$ be a domain. The configuration space $X\subset L^2(\Omega,\Cb)$ is taken
		to be a suitable space of wave functions,
		and we equip $X$ with the real inner product $\inn{\psi}{\phi}=\int\Re\b\psi\phi$.
		Then the map $\psi\mapsto (\Re \psi,\Im \psi)$ is an isometric isomorphism 
		between $X$ and a subspace of  $L^2(\Omega,\Rb)\times  L^2(\Omega,\Rb)$,
		if the latter is equipped with inner product
		$\inn{u}{v}=\tfrac{1}{2}\int (u_1v_1+u_2v_2).$
		Hence, we can identify $X$ as a subspace of $L^2(\Omega,\Rb)\times  L^2(\Omega,\Rb)$.
		
		Under this identification, the operator 
		$J:\psi\mapsto -i\psi$ can be represented by the symplectic matrix
		\begin{equation}
			\label{4.0}
			J=\begin{pmatrix}0&1\\-1&0\end{pmatrix}.
		\end{equation}
		Consider now a Schr\"odinger type equation
		\begin{equation}
			\label{4.1}
			i\di_t u=E'(u),
		\end{equation}
		where, for consistency of notation, $E$ is a suitable Hamiltonian.
		Using the identification described above, we can cast \eqref{4.1}
		into \eqref{Ham}.
		This is the typical setting for evolutions from quantum mechanics. 
		
		2. Let $X_1\subset L^2(\Omega,\Rb^k)$ be 
		a suitable space of density functions or order parameters,
		equipped with the usual inner product. 
		
		Consider a second-order dynamics
		\begin{equation}
			\label{4.2}\di_{tt}v=-\tilde{E}'(v),\end{equation}
		where $v\in X_1$ and $\tilde{E}$ is some Hamiltonian.
		We can reduce this to a first-order system 
		by setting $u=(u_1,u_2)=(v,\di_tv)$,
		and choosing some $E$ s.th. 
		$E'(u)=(\tilde{E}'(u_1), u_2)$.
		The configuration space for $u$ is $X=X_1\times X_2$,
		for some suitable $X_2\subset L^2(\Omega,\Rb^k)$.
		$X$ is equipped with the inner product 
		$\inn{\cdot}{\cdot}_X=\inn{\cdot}{\cdot}_{X_1}+\inn{\cdot}{\cdot}_{X_2}.$
		This way we can cast \eqref{4.2} into \eqref{Ham},
		with $J$ given by \eqref{4.0}.
		
		The main result of this section is the next theorem,
		analogous to \thmref{thm1} in the gradient case, but only valid on a 
		long finite time interval.
		
		\begin{theorem}
			\label{thm2}
			Fix any $T\gg1$ and $0<\g<1$. There exists $0<\eps_0\ll1$ depending on $\g,\,T$
			only, s.th. the following holds:
			For every $0<\eps\le\eps_0$,  there exists $0<\delta\ll\eps$
			s.th.
			the following holds:
			Let $M$ be the manifold of approximate solitons
			as in \defnref{M}.
			Let $u_0\in X$ be an initial configuration s.th.
			$\dist(u_0,M)\le\delta$. Let $u_t$ be the flow
			generated by $u_0$ under \eqref{Ham}.
			\begin{enumerate}
				\item (A priori estimate) For all $\epsilon^\g t\le T$,  there holds
			\begin{equation}
				\label{4.3}
				\dist(u_t,M)\le\eps^\g.
			\end{equation}
			
			\item (Effective dynamics) Moreover, for $\eps^\g t\le T$,
			there exists 
			a unique decomposition for $u_t$ as in \eqref{2.5.3}, 
			and the moduli $\s\equiv \si_t:=S(u_t)$ satisfies the following effective
			dynamics:
			\begin{equation}
				\label{4.4}
				\di_t\s=\J_\s^{-1}\E'(\s)+O_{\norm{\cdot}_Y}(\eps^{1+\al+\g}).
			\end{equation}
			\end{enumerate}

		\end{theorem}
		\begin{remark}
			As in the gradient case, the remainder in \eqref{4.4} is of lower order by
			\eqref{5.8},
			but here we have a weaker error estimate.
			Note importantly that the implicit constant in the remainder of \eqref{4.4} is independent of $T$ and $\g$.
			
		\end{remark}
		\begin{proof}[Proof of \thmref{thm2}]

			1. The initial setup is identical to the first step in \thmref{thm1}.
			In particular, we write $u_t=v_t+w_t$ as in \eqref{decomp}. 
			
			The claim now is that for every $ t\le \eps^{-\g}T$,
			there holds
			\begin{equation}\label{A4.1}
				\norm{w_t}_X\le \eps^\g,
			\end{equation}
			provided $0<\delta\ll\eps\le\eps_0(T,\g)$.

			Notice that 
			by the continuity of \eqref{Ham},
			there exists some
			(possibly small) $T_1>0$ s.th. \eqref{A4.1} holds
			for all $t\le T_1$,
			provided $\norm{w_0}_X\ll \eps^\g$.
			The latter is the case because $0<\g<1$, and for $\delta\ll\eps$, we have 
			$\norm{w_0}_X=O(\eps)$ by \remref{remD}. 
			Hence, the decomposition $u_t=v_t+w_t$
			is valid at least locally.
			
			As in Step 1 of \thmref{thm1},
			if \eqref{A4.1} holds, then \eqref{4.3} follows
			 from definition \eqref{dist}.
			By \lemref{lem2.1}, \eqref{4.3} implies the validity of
			the decomposition \eqref{2.5.3} for $u_t$.
			By  Lemma \ref{lem2.2}, \eqref{A4.1} implies 
			the effective dynamics \eqref{4.4} and the remainder estimate therein.

			Hence, the theorem
			is proved once we establish \eqref{A4.1} for all $ t\le \eps^{-\g}T$.
			
			2. Consider the expansion
			\begin{equation}
				\label{3.12}
				E(u)=	E(v+w)=E(v)+\inn{E'(v)}{w}+\frac{1}{2}\inn{L_\s w}{w}+R_\s(w),
			\end{equation}
			where $R_\s(w)$ is the super-quadratic remainder.
			By the construction from \lemref{lem2.1}, the fluctuation field
			$w$ satisfies \eqref{2.5}. Thus $w\in (JT_vM)^\perp$, and 
			by condition \eqref{C2}, we can rearrange \eqref{3.12} to obtain
			\begin{equation}
				\label{3.12.1}
				\norm{w}_{X}^2\ls E(v+w)-E(v)-\inn{E'(v)}{w}-R_\s(w).
			\end{equation}
			
			Since $E(u)$ is conserved along \eqref{Ham}, we have
			$$E(v+w)=E(v_0+w_0)=E(v_0)+\inn{E'(v_0)}{w_0}+\frac{1}{2}\inn{L_\s w_0}{w_0}+R_{\s_0}(w_0).$$
			Plugging this into \eqref{3.12.1}, we have
			\begin{equation}
				\label{3.12.2}
				\begin{aligned}
					\norm{w}_{X}^2\ls& E(v_0)-E(v)\\&+\inn{E'(v_0)}{w_0}-\inn{E'(v)}{w}+\frac{1}{2}\inn{L_{\s_0} w_0}{w_0}+R_{\s_0}(w_0)-R_\s(w).
				\end{aligned}
			\end{equation}
			
			The last five terms in \eqref{3.12.2} can be controlled as follows:
			By the approximate critical point property \eqref{C1},
			we have $$\inn{E'(v_0)}{w_0}-\inn{E'(v)}{w}\ls \eps(\norm{w_0}_X+\norm{w}_X)\ls \eps M(t).$$
			Here recall that we have defined the function
			$$M(t):=\sup_{t'\le t}\norm{w_{t'}}_{X}.$$
			
			Next, since $L_{\s_0}$ is bounded,
			we have 	
			$$\inn{L_{\s_0} w_0}{w_0}\ls\norm{w_0}_X^2.$$
			Note that this bound does not depend on $\si_0$ (see the discussion about \eqref{20'}).
			
			Lastly, since the remainder $R_\s(w)$  is of the order $o(\norm{w}_X^2)$ for $C^2$ functional $E$, it follows that 
			$$R_\s(w)-R_{\s_0}(w_0)=o(M(t)^2).$$
			By the preceding estimates,  \eqref{3.12.2} becomes
			\begin{equation}
				\label{3.12.3}
				\norm{w}_{X}^2\ls E(v_0)-E(v)+
				\eps M(t)
				+o(M(t)^2)+\norm{w_0}_X^2.
			\end{equation}
			
			3. It remains to control the first two terms in the r.h.s. of \eqref{3.12.3}.
			This difference is the energy fluctuation of the approximate solitons,
			and it can be controlled as follows.
			Differentiating the energy $E(t)=E(f(\s_t))$ and using \eqref{3.7}, we have
			\begin{equation}
				\label{3.13}
				\begin{aligned}
					\od{E}{t}&=\inn{E'(v)}{\di_tv}\\
					&=\inn{E'(v)}{Q_\s JE'(v)}+\inn{E'(v)}{Q_\s (JL_\s w-\di_tw)}+\inn{E'(v)}{Q_\s JN_\s(w))}.
				\end{aligned}
			\end{equation}
			We now bound the three inner products in the last line. 
			
			Recall the definition of $Q_\s$  in \eqref{2.4} as the skew projection, which gives
			the relations  $Q_\s^2=Q_\s$ and $Q_\s J=JQ_\s^*$ (the latter follows from \eqref{J} and \eqref{2.4'}).
			Using these 
			and the fact that $J$ is symplectic, we find
			$$\begin{aligned}
				\inn{\phi}{Q_\sigma J\phi}&=\inn{\phi}{Q_\sigma^2 J\phi}\\
				&=\inn{Q_\sigma ^*\phi}{Q_\sigma J\phi}\\
				&=\inn{(J^{-1}J)Q_\sigma^* \phi}{Q_\sigma J \phi}\\
				&=-\inn{J(Q_\sigma J \phi)}{Q_\sigma J \phi}=0\quad (\phi\in X).
			\end{aligned}$$
			Applying this with $\phi=E'(v)$, we see that  the first term in \eqref{3.13} vanishes. 
			
			For the two estimates below, we need to use the approximate critical point
			property \eqref{C1}.
			Recall also that as discussed in \remref{remD},
			we can drop higher order terms in $w$
			on the interval $\eps^\g t\le T'$,
			so long as \eqref{A4.1} holds.

			Using \eqref{3.9}-\eqref{3.10} and \eqref{2.8}, the second inner product in the last line of  \eqref{3.13} can be bounded as
			\begin{equation}\label{3.14}\abs{\inn{E'(v)}{Q_\s(JL_\s w-\di_tw)}}\ls\eps^2\norm{w}_X.
			\end{equation}

			By the uniform nonlinear estimate \eqref{3.11},
			the third inner product  can be bounded as
			\begin{equation}
				\label{3.15}
				\abs{\inn{E'(v)}{Q_\s JN_\s(w))}}\ls\eps\norm{w}_X^2.
			\end{equation}
			Combining \eqref{3.13}-\eqref{3.15} and integrating from $0$ to $t$, we find
			\begin{equation}
				\label{3.16}
				\abs{E(v(t))-E(v(0))}\ls t\del{\eps^2M(t)+ \eps M(t)^2}.
			\end{equation}
			
			4. Plugging \eqref{3.16} into
			\eqref{3.12.3}, and then dividing both side by $M(t)$, we have 
			\begin{equation}
				\label{3.17}
				M(t)\ls t\del{\eps^2+ \eps M(t)} + \eps +o(M(t))+\norm{w_0}_X.
			\end{equation}
			For $\delta\ll\epsilon$, the last term is  $O(\epsilon)$ by \remref{remD}.
			Hence, there exist two constants 
			$C_0>0$ and $0<\delta_0\ll1$,
			both independent of $t,\,T,\,\g$,
			such that
\begin{equation}\label{3.17'}
				M(t)\le C_0\del{t(\eps^2+\eps M(t)+\eps+\delta_0M(t))+\eps}.
\end{equation}
			Now, choose a sufficiently small $\eps_0=\eps_0(\g,T)>0$ from here,
			so that 
			$$1-C_0(t\eps+\delta_0t)>1/2\quad (t\le \eps^{-\g}T).$$
			Then it follows from \eqref{3.17'} that 
\begin{equation}\label{3.17''}
				M(t)\le 2C_0T(\eps^{2-\g}+\eps)+2C_0\eps\quad (t\le \eps^{-\g}T).
\end{equation}
%

			At this point, the leading order term in
			the r.h.s. of \eqref{3.17''} is of the order $\eps$.	
			Since $0<\g<1$, and $C_0$ is an absolute constant, we conclude from \eqref{3.17''} that 
			we can further shrink $\eps_0=\eps_0(\g,\,T)>0$ so that  
			\eqref{A4.1} holds for every  $\eps^\g t\le T$ with $0<\eps\le\eps_0$.
			This completes the proof.

		\end{proof}

		\begin{corollary}[Converse of Theorem \ref{thm2}]\label{cor2}
			Fix any large $T\gg1$ and $0<\g<1$.
			There exists $c>0$ independent of $T$ and 
			$0<\eps_0\ll1$ depending on $\g,\,T$ only 
			with the following properties:
			Let $0<\eps\le\eps_0$.
			Let $\s_0\in\Sigma$. Let $\s_t\in\Sigma$ be the flow generated by $\s_0$
			under the effective dynamics \eqref{Ham}. Then there exists a solution $u_t$ to \eqref{3.1}
			s.th. for all $\epsilon^\g t\le T$, there holds
			\begin{equation}
				\label{23'}
				\norm{u_t-f(\sigma_t)}_X\le c \max(\epsilon^{(1-\g)/2},\eps^\g).
			\end{equation}
		\end{corollary}
		\begin{proof}
			Choose $\eps_1=\eps_1(\g,T)>0$
			s.th. \thmref{thm2} holds with $\g,\,T$ and all $0<\eps\le\eps_1$.
			Following the proof of Corollary \ref{cor1}, \textit{mutatis mutandis},
			we find 
			\begin{equation}\label{12121}
				\norm{u_t-f(\sigma_t)}_X\le c_1T(\eps^{1-\g}+c_2\eps) + c_3\eps^\g \quad (\eps^\g t\le T).
			\end{equation}
			Here $c_1,\,c_2,\,c_3>0$ are all absolute constants. 
			Since $\g<1,$ we can choose $0<\eps_0(\g,\,T)\le \eps_1$ from \eqref{12121}
			to conclude \eqref{23'}.
		\end{proof}
		
		\section{Application to Interface dynamics}
		\label{sec:5}
		In this section we consider a typical situation arising
		from the study of phase transition.
		Following \cite{VortFil}, we give some concrete examples of approximate solitons
		that fully utilize the generality of the adiabatic framework developed in the 
		preceding sections.
		
		Consider a suitable function space $X$ consisting of vector-valued order parameters $$\psi:\Rb^n_x\times\Rb^k_z\to \Rb^k\quad (n\ge1,\,k\ge0).$$
		Suppose there is a $k$-dimensional interface at $\Set{x=0,z\in\Rb^k}$
		separating two homogeneous phases with a steep transition layer. 
		Here $k=0$ is allowed, because one can be  interested
		in some soliton concentrated
		at a single point in $\Rb^n$ (e.g. one-dimensional kinks,
		planar Ginzburg-Landau vortices, spherical droplets around some point in the space,
		etc.).
		
		Suppose
		we are given a stable equilibrium $\psi$ of \eqref{1}, s.th. for some $0\le\alpha<1$ and small $\eps\ll1$,  
		\begin{align}
			\label{5.3}
			&\norm{\grad_x\psi}_X\sim \eps^{-\alpha}\text{ due to the steep phase transition at the interface},\\
			\label{5.9}
			\begin{split}
				&E'(\psi)=0,\quad L_\psi:=E''(\psi)\ge0,\\
				&\text{$0$ is isolated from the rest of the spectrum of $L_\psi$}. 
			\end{split}
		\end{align}
		Note that we do not require $L_\psi>0$, because if there are 
		continuous symmetries broken by $\psi$, then
		$L_\psi$ in general have zero modes due to symmetry breaking.

		Consider another  space $Y$ consisting of
		smooth perturbations of the form
		\begin{equation}\label{si}
			\s:\Rb^k_z\to \Rb^n,\quad \norm\si_Y\ll1.
		\end{equation}
		Geometrically, such $\si$ can be thought of as some  ``wiggling'' within the $n$-dimensional horizontal
		cross sections around
		the $k$-dimensional interface $\Set{x=0}$.
		
		Define  a map
		\begin{equation}
			\label{5.1}
			f:\s\mapsto \psi_\s:=\psi(x-\s(z),z).
		\end{equation}
		This map is smooth if $\psi$ is smooth. 
		Below we consider the two groups of assumptions from \secref{sec:2}
		in connection with this particular parametrization \eqref{5.1}.
		
		For the first group of assumptions, we note that 
		by \eqref{5.9}, for $\si=0$, the configuration
		$f(0)=\psi$ is an exact soliton satisfying conditions \eqref{C1}-\eqref{C3}
		with $\eps=0$. 
		As we discussed in \secref{sec:2.1},
		these conditions are persistent for small
		perturbations of a given approximate soliton,
		and it follows that 
		\eqref{C1}-\eqref{C3} hold
		on an appropriately chosen manifold $\Si$ around $\si=0$ consisting
		of perturbations of the form \eqref{si}.
		The choice  of $\Si$ depends only on
		$\ker L_\psi$ for the linearized operator at the given 
		exact soliton in \eqref{5.9}.
		See e.g. \cite{VortFil} in  the exactly same setting, and
		\cites{MR1257242,MR1309545,MR2097576,MR2189216} in closely related settings.

		For the second group of geometric assumptions, we note that  
		the Fr\'echet derivative of $f$ is given by 
		\begin{equation}
			\label{5.2}
			df(\s):\xi \mapsto -\grad_x\psi_\s\cdot \xi\quad (\xi\in Y).
		\end{equation}
		This map is clearly injective, so $f$ is an immersion and \eqref{5.4'} holds.
		Using the formula \eqref{5.2},
		depending on the particular choice of 
		spaces $X$ and $Y$, 
		one can deduce \eqref{5.4} from  the condition  \eqref{5.3}.
		
		To illustrate this fact, consider the situation $n=1,\, k\ge 1$.
		Fix  the configuration spaces
		$$X:=L^2(\Rb_x\times \Rb^k_z,\,\Rb^k),\quad Y:=L^2(\Rb^k_z,\Rb).$$
		Suppose we are given a function $\tilde\psi:\Rb_x\to \Rb^k$ satisfying 
		\begin{equation}\label{6.1}
			\int_\Rb \abs{\grad_x \tilde\psi(x)}^2\,dx=A^2,
		\end{equation}
		e.g. a certain curve of geometric interest. 
		Then we define the  lift $\psi:\Rb_x\times \Rb_z^k\to \Rb^k$ 
		by setting $\psi(x,z)=\tilde\psi(x)$ for every $z$.
		In principle, the solitonic properties 
		of $\tilde\psi$ is not affected by such a lift,
		see e.g. \cite{VortFil}.
		
		Suppose the parametrization  $f$ and its derivative is
		given by \eqref{5.1}-\eqref{5.2}, and recall
		$g_\si:Y\to X$ is defined as the 
		Fr\'echet derivative $df(\si)$ in local coordinate.
		For $\xi\in Y$, we compute
		\begin{equation}\label{6.2}
			\begin{aligned}
				\norm{g_\si\xi}_X^2&=\int_{\Rb}\int_{\Rb^k}\abs{\xi(z)\cdot\grad_x\psi_\si(x,z)}^2\,dxdz\\
				&=\int_{\Rb^k}\abs{\xi(z)}^2\del{\int_\Rb \abs{\grad_x\psi(x-\s(z),z)}^2\,dx}\,dz\\
				&=\int_{\Rb^k}\abs{\xi(z)}^2\del{\int_\Rb \abs{\grad_x\tilde\psi(x-\s(z))}^2\,dx}\,dz\\
				&=\int_{\Rb^k} \abs{\xi(z)}^2A^2\,dz=\norm{\xi}_Y^2A^2.
			\end{aligned}
		\end{equation}
		In the second equality we use Fubini's theorem, and in the fourth equality we use
		the fact that  Lebesgue measure is translation invariant. 
		
		Now, if $A\sim \eps^{-\al}$ in \eqref{6.1}, then we can conclude
		\eqref{5.4} from \eqref{6.2}.
		Thus, with this example,
		we have demonstrated the typical implication from \eqref{5.3} to \eqref{5.4}.
		
		\begin{remark}
			This kind of geometric implication 
			was already noted and played a crucial role in \cite{MR2189216} (see Sect. 3.1
			of that paper), as well as \cites{MR1257242,MR1309545},
			where the authors studied  some situations with $n=2,\, k=0$.
		\end{remark}
		
		The conclusion from the above discussion	is that  
		the framework we laid out in the previous sections applies
		to the study of interface dynamics via parametrization of the form \eqref{5.1},
		provided there is some known exact soliton to begin with.
		We proceed to demonstrate in the next subsection an  
		application with $n=2,\,k=1$.
		
		\subsection{Example: Ginzburg-Landau Vortex Filaments}\label{sec:6.1}
		In this subsection, we discuss some of the results obtained in \cite{VortFil}
		in the adiabatic framework developed in the present paper.
		
		Let $\Omega\subset \Rb^d,\,d=2,3$ be a  domain.
		Consider the Ginzburg-Landau energy functional 
		\begin{equation}
			\label{GL}
			E(\psi)\equiv E_\Omega^\eps(\psi):=\int_\Omega \frac{1}{2}\abs{\grad\psi}^2+\frac{1}{4\eps^2}\del{\abs{\psi}^2-1}^2.
		\end{equation}
		Here $\psi:\Omega\to\Cb$ is a complex order parameter representing,
		for instance, the Bose-Einstein 
		condensate in superfluidity. 
		The energy \eqref{GL} has translation, rotatoin, and global  $U(1)$-gauge symmetries.
		
		It is well-known that there exist non-trivial stable critical points
		$\tilde\psi:\Omega\subset \Rb^2\to\Cb$ satisfying 
		\begin{align}
			&\label{5.10}
			\norm{\grad_x\tilde\psi}_{L^2}\sim \abs{\log\eps}^{1/2},\\
			&\label{5.11}
			E'(\tilde\psi)=0,\quad L_{\tilde\psi}\ge0,\quad {L_{\tilde\psi}}\vert_{Z^\perp}\ge\beta>0.
		\end{align}
		See for instance \cites{MR1269538,MR1763040}. 
		Here $Z$ denotes the space of symmetry zero modes, in this case generated by
		the broken translation and global gauge symmetry. Such $\tilde\psi$ are 
		known as the (planar) \emph{vortex} solutions.
		
		The characteristic feature of a vortex $\tilde\psi$ is its concentration 
		property.
		This is due to the structure of nonlinearity  in \eqref{GL}.
		The quartic, hat-shaped potential term 
		forces the modulus $\abs{\tilde\psi}$ of the non-homogeneous equilibrium to rapidly increase from $0$ to 
		$1$ in all directions away from 
		the vortex center.
		
		The planar vortex configuration $\tilde\psi$ obviously lifts to a steady state  in $\Rb^3$ through
		$\tilde\psi(x)\mapsto \psi(x,z)\equiv\tilde\psi(x)$, where $z\in\Rb$ parametrizes 
		the vertical direction. 
		This lift $\psi(x,z)$
		concentrates near the vertical axis $\Set{x=0}$. 
		
		Now, consider 
		perturbations of the form (up to a global gauge)
		\begin{equation}
			\label{5.12}
			f(\s )=\psi(x-\s(z),z),\quad \s:\Rb\to\Rb^2.
		\end{equation}
		Since each planar vortex $\tilde\psi(x)$ concentrates around $x=0$,
		the function $f(\s)$ describes a vortex filament that ``curves around'' a 
		concentration set near 
		the 
		vertical axis $\Set{(0,z)\in\Rb^3}$. 
		In general, \eqref{5.12} is not a critical point of \eqref{GL},
		and do not arise from any symmetry reduction procedure.
		Moreover, the space of $\si$ can be infinite-dimensional.
		
		In \cite{VortFil}, we
		show that under some small curvature
		assumption on the perturbation parameter $\si$
		(which lifts to a three-dimensional curve through $\si(z)\mapsto(\si(z),z)$ that winds 
		around the axis $\Set{(0,z)\in\Rb^3}$),
		the parametrization \eqref{5.12} (up to a global gauge)
		gives a manifold of approximate solitons 
		sitting in the energy space for \eqref{GL},  satisfying all the assumptions in \secref{sec:2}.
		
		Using the method developed in \secref{sec:H} for Hamiltonian system,
		we obtain an adiabatic approximation for  the three-dimensional Gross-Pitaevskii equation 
		\begin{equation}\label{5.13}
			i\pd{\psi}{t}=-\Lap \psi+\frac{1}{\eps^2}\del{\abs{\psi}^2-1}\psi\quad (\psi:\Omega\subset \Rb^3\to\Cb).
		\end{equation}
		Evolution \eqref{5.13}
		can be cast into a Hamiltonian system of the form \eqref{4.1}, with $E$ given
		by the Ginzburg-Landau energy \eqref{GL}.
		The effective dynamics for \eqref{5.13} is an evolution of the concentration sets $\s$ (i.e. the moduli, using the terminology in \secref{sec:2}), namely the binormal 
		curvature flow
		\begin{equation}\label{BNF}
			\di_t\vec\s=\di_s\vec\s\times \di_{ss
			}\vec\s,
		\end{equation}
		where \( \vec\si(z):=(\si(z),z)\) is the lift of $\si$ to a spatial curve, and $s=s(z)$ is the arclength parameter, satisfying $\tfrac{ds}{dz}=\abs{\vec\si}$.
		The flow \eqref{BNF} is a Hamiltonian system in 
		the moduli space (in this case consists of functions  $\si:\Rb\to \Rb^2$), and it appears in place of the abstract effective dynamics \eqref{eff}, \eqref{4.4}.
		
		Let us remark that most  results relating the geometry of $k$-dimensional 
		interfaces with $k\ge1$  to the full
		configurations on $\Rb^{n+k}$ rely on rather involved measure
		theoretic arguments. See for instance the important 
		contributions \cites{MR1491752,MR2195132}. 
		Moreover, these results do not retain the structure
		of the interface, as they take the limit as the length scale 
		$\eps\to 0$.
		On the other hand, results using adiabatic approximations are mostly
		for rather simple geometry of interfaces (e.g. 
		finite collection of points, in which case $k=0$), see the cited works in  
		Introduction. 
		
		\section*{Acknowledgment}
		The Author is supported by Danish National Research Foundation grant CPH-GEOTOP-DNRF151.  
		
		\section*{Declarations}
		\begin{itemize}
			\item Conflict of interest: The Author has no conflicts of interest to declare that are relevant to the content of this article.
			\item Data availability: Data sharing is not applicable to this article as no datasets were generated or analysed during the current study.
		\end{itemize}

		\appendix
		\section{Basic Variational Calculus}
		Here we recall some basic elements of variational calculus that have been used 
		repeatedly. For details, see for instance \cite{MR2431434}*{Appendix C}, \cite{MR1336591}*{Chapt. 1}.
		
		\subsection{Fr\'echet Derivative} Let $X,\,Y$ be two Banach spaces. Let $U$ be an open set in $X$.
		For a map $g:U\subset X\to Y$ and a vector $u\in U$,
		the Fr\'echet derivative $dg(u)$ is a linear map from $X\to Y$ s.th.
		$g(u+v)-g(u)-dg(u)v=o(\norm{v}_X)$ for every $v\in X$ with $\norm{v}_X\ll1$ . 
		If $dg(u)$ exists at $u$, then it is unique. If $dg(u)$ exists for every $u\in U$,
		and the map $u\mapsto dg(u)$ is continuous from $U$ to the space of linear operators $L(X,Y)$, then we
		we say $g$ is $C^1$ on $U$. 
		In this case, $dg(u)$ is uniquely given by 
		$$v\mapsto {\od{g(u+tv)}{t}}\vert_{t=0}\quad (v\in X).$$
		Iteratively, we can define higher order derivatives this way.

		\subsection{Gradient and Hessian} If $X$ is a Hilbert space over a scalar field $Y$, then by Riesz representation,
		we can identify $dg(u)$ as an element in $X$, denoted by $g'(u)$. The vector
		$g'(u)$ is called the $X$-gradient of $g$. 
		Similarly, we denote $g''(u)$ the second-order Fr\'echet derivative
		$d^2g(u)$. If $g$ is $C^2$, then $g''$ can be identified as a symmetric 
		linear operator
		uniquely determined  by the relation  
		$$\inn{g''(u)v}{w}={\md{g(u+tv+sw)}{2}{t}{}{s}{}}\vert_{s=t=0}\quad (v,w\in X).$$
		
		\subsection{Remainder and Composition} Let $X$ be a Hilbert space over a scalar field $Y$. Suppose $g$ is $C^2$ on $U\subset X$. Define a scalar function $\phi(t):=g(v+tw)$ for vectors
		$v,w$ s.th. $v+tw\in U$ for every $0\le t\le 1$.
		Then the elementary Taylor expansion at $\phi(1)$ gives 
		$$g(v+w)=g(v)+\inn{g'(v)}{w}+\frac{1}{2}\inn{g''(v)w}{w}+o(\norm{w}_X^2).$$
		Here we have used the definition of $g'$ and $g''$ from the last subsection.
		
		Let $\Omega\subset \Rb^d$ be a bounded domain with smooth boundary.
		Fix $r>d/2,\,f\in C^{r+1}(\Rb^n)$. For $u:\Omega\to \Rb^n$,
		define a map $g:u\mapsto f\circ u$. Then $g:H^r(\Omega)\to H^r(\Omega)$
		is $C^1$, 
		and the Fr\'echet derivative is given by $v\mapsto \grad f\cdot v$.

		\bibliography{bibfile}
	\end{document}